\documentclass[12pt]{article}
\usepackage[english]{babel}
\usepackage[latin1]{inputenc}
\usepackage{amsfonts,amssymb,amsmath, epsfig}
\usepackage{color,graphicx,graphics,psfrag}
\usepackage{amsmath,amstext,amssymb,amsfonts, amscd}
\usepackage{hyperref}           
\usepackage{mathtools}   
\usepackage{authblk}

\newcommand{\RR}{\mathbb{R}}

\definecolor{cinnamon}{rgb}{0.82, 0.41, 0.12}
\def\Box{\leavevmode\vbox{\hrule
     \hbox{\vrule\kern4pt\vbox{\kern4pt}%
           \vrule}\hrule}}
\def\blackbox{\leavevmode\vrule height 5pt width 4pt depth 0pt\relax}
  
\def\endprf{\null \hfill {$\blackbox$} \bigskip}

\newcounter{appendix}
\def\appendix{\advance\c@appendix by 1
   \def\thesection{\Alph{section}}
   \ifnum\c@appendix=1 \setcounter{section}{-1} \fi
   \@startsection {section}{1}{\z@}{-3.5ex plus -1ex minus 
   -.2ex}{2.3ex plus .2ex}{\Large\bf}}

\def\paragraph#1{{\bf #1\ }}

\newtheorem{lemma}{Lemma}[section]  

\newtheorem{theorem}[lemma]{Theorem}
\newtheorem{proof}[lemma]{Proof}

\newtheorem{definition}[lemma]{Definition}

\newtheorem{proposition}[lemma]{Proposition}

\newtheorem{remark}{Remark}[section]

\newcommand{\R}{\mathbb{R}}

\title{A new continuum theory for incompressible \\swelling materials}

\author[(1)]{Pierre Degond}
\author[(2)]{Marina A. Ferreira}
\author[(3)]{Sara Merino-Aceituno}
\author[(4)]{ Micka\"el Nahon}

\affil[(1,2,3)]{Department of Mathematics, Imperial College London, South Kensington Campus,
London, SW7 2AZ,
UK}
\affil[(1)]{pdegond@imperial.ac.uk}
\affil[(2)]{m.amado-ferreira14@imperial.ac.uk}
\affil[(3)]{s.merino-aceituno@imperial.ac.uk}
\affil[(4)]{\'Ecole Normale Sup\'erieure de Lyon, 46 All\'ee d'Italie, 69007 Lyon, France\\
mickael.nahon@ens-lyon.fr}
 
\date{}

\begin{document}

\maketitle


\vspace{0.5 cm}
\begin{abstract}
Swelling media (e.g. gels, tumors) are usually described by mechanical constitutive laws (e.g. Hooke or Darcy laws). However, constitutive relations of real swelling media are not well-known. Here, we take an opposite route and consider a simple packing heuristics, i.e.  the particles can't overlap. We deduce a formula for the equilibrium density under a confining potential.  We then consider its evolution when the average particle volume and confining potential depend on time under two additional heuristics: (i) any two particles can't swap their position; (ii) motion should obey some energy minimization principle. These heuristics determine the medium velocity consistently with the continuity equation. In the direction normal to the potential level sets the velocity is related with that of the level sets while in the parallel direction, it is determined by a Laplace-Beltrami operator on these sets. This complex geometrical feature cannot be recovered using a simple Darcy law.

\end{abstract}

\medskip
\noindent
{\bf Acknowledgements:} PD acknowledges support by the Engineering and Physical Sciences 
Research Council (EPSRC) under grants no. EP/M006883/1 and EP/N014529/1, by the Royal 
Society and the Wolfson Foundation through a Royal Society Wolfson  Research Merit Award no. WM130048 and by the National Science  Foundation (NSF) under grant no. RNMS11-07444 (KI-Net). PD is on leave from CNRS, Institut de Math\'ematiques de Toulouse, France.\\
S.M.A. was supported by the British Engineering and Physical Research Council under grant ref: EP/M006883/1.\\
MF acknowledges support of the Department of Mathematics, Imperial College London, through a Roth PhD scholarship. \\
 MN gratefully acknowledges the hospitality of the Department of Mathematics, Imperial College London, where this research was conducted. 

\medskip
\noindent
{\bf Data statement: } No new data were collected in the course of this research. 

\medskip
\noindent
{\bf Conflict of interest:} The authors declare that they have no conflict of interest. 


\vskip 0.4cm

\setcounter{equation}{0}
\section{Introduction}
\label{sec_intro}

Swelling or drying media are encountered in many contexts such as chemistry or material science (swelling gels), biology (cancer tumors or growing tissues), geosciences (drying of wetting soil), cooking (dough being cooked), etc. The modelling of swelling or drying media from first principles is difficult due to the complex nature of the materials (cells, mixtures, polymers, etc). Often, they have intermediate properties between solids and liquids or can have genuinely new properties (biological tissues). Modelling of swelling or drying material is very important in view of potential applications in health (tumour growth or tissue development) and other sciences. 

Modelling of swelling material can be attempted through either solid or fluid mechanics models. In the first category, we refer to \cite{amar2010swelling} (and references therein) where a model of swelling gel is proposed in the framework of hyperelasticity theory. Interestingly, this model was developed as swelling gels are seen as a good laboratory model of certain tumors, such as malignant melanoma. Indeed, the instabilities that are observed at the boundary of the gel are reminiscent to the corrugated shape of the boundary of a melanoma. In the context of tumor growth modelling a solid mechanics models can be found e.g. in \cite{chaplain1992mathematical}. 

However, many of the models used in tumor growth rather use a fluid-dynamic approach, and specifically, Darcy's law or some elaboration of it \cite{bertsch2012nonlinear, byrne2003modelling, cristini2003nonlinear, hawkins2012numerical, colli2015vanishing}.  Mathematically, Darcy's law is expressed by $v= -k\, \nabla p$, where $v$ is the fluid velocity, $p$ is the hydrostatic pressure, $\nabla$ is the spatial gradient and $k$ is a constant named 'hydraulic conductivity'. Darcy's law is derived from Navier-Stokes equation for a fluid subjected to strong friction such as flowing inside a porous medium. However, the use of Darcy's law is not obvious. The article \cite{ambrosi2002closure} is entirely devoted to the problem of determining the velocity in the mass balance equations (referred to as the ``closure problem'') and to a phenomenological justification of the use of Darcy's law in tumour growth. 

Due to its importance in the clinic, one of the major questions explored in tumor growth modelling is the description of the tumor boundary and how it evolves in time. It naturally leads to the study of free boundary problems \cite{friedman2006free} and many works have explored under which asymptotic limits the fluid model could lead to a free-boundary problem \cite{bertsch2010free, hilhorst2015formal}. Related to these, the analogy between tumor growth and the free-boundary problem of solidification (the so-called Hele-Shaw problem) has been developed in \cite{hecht2017incompressible, perthame2014hele, perthame2014derivation, perthame2015incompressible}. In these last series of works, the tumor is regarded as the region of space where cells have reached the packing density. It presupposes that the cells have a finite size and cannot overlap, leading to a maximum packing density where cells occupy all the available space. The tumor is therefore an incompressible medium separated from the outer medium by a moving free boundary which can be calculated through the resolution of an elliptic problem for the pressure in the moving domain of the tumor. 

All the previous studies rely on a continuum description of the tumor. However, at the microscopic level, a tumor is made of discrete entities, the cells and various types of ``individual-based'' microscopic models of tumor growth, where cells are described as discrete entities, have been developed: see in particular \cite{drasdo2005single}. We refer to \cite{roose2007mathematical} for a review of the various modelling approaches and to \cite{byrne2009individual} for a comparison of their merits. The connection of the microscopic approach to the macroscopic one through coarse-graining is investigated in \cite{motsch2017short}. 

In the present work, we revisit the closure problem and investigate what motion results from the combination of volume-exclusion (or non-overlapping) and growth. In relation to this, we question the validity of Darcy's law once more. Our approach, rather than relying on constitutive relations like hyper-elasticity or Darcy's law, hypothesizes simple heuristic rules, more likely to be obeyed in generic situations. Here, the main heuristic rule is that particles cannot overlap. In other words, we directly place ourselves in a context akin to the Hele-Shaw limit as developed in \cite{perthame2014hele} and related works cited above. However, as we will see, our conclusions will be different. We also point out that similar heuristic rules have been applied to other domains, such as crowd modelling (see in particular \cite{roudneff2011handling}). 

We consider a system made of finite-sized particles at equilibrium in a confining external potential constrained by the non-overlapping condition. We refer to \cite{leroy2017tumor} for a discussion of the biological situation described by this particular setting. We then assume that the particle volume and confinement potential may vary with time and that the particles follow this evolution adiabatically by remaining at any time at mechanical equilibrium. The question we want to address is what particle motion results from this situation. 

Answering this question in full generality at the discrete level is probably out of reach. So, we formulate a similar problem at the continuum level. We assume a continuum density for a population of particles having finite average volume. The particles are confined by an external potential and we assume the particles at mechanical equilibrium. Our first result is to characterize the resulting equilibrium density. Like in the Hele-Shaw type models referred above, the particles occupy a domain of finite extension in space, limited by a level set of the potential. Inside this domain, the density is equal to the maximal (packing) density allowed by their finite size. Outside this domain the density is zero. 

Then, we turn on the time variability of the average particle volume and of the confinement potential. Assuming that the system moves adiabatically and remains at any time at mechanical equilibrium, we can compute the continuum velocity. More precisely, we determine this velocity by applying two heuristic principles directly connected to the previous non-overlapping heuristics. The first heuristics is that particles can't swap their positions. Indeed, at the packing state, there is not enough space for two spherical particles to undertake the maneuver required to swap their position. This heuristics provides the component of the velocity normal to the potential level sets.

To determine the component of the velocity tangent to the potential level sets, we invoke a second heuristics, namely that the sequence of minimization problems over time will favor a continuous particle motion rather than jumps which would generate large velocities. In continuum language, this means that the velocity should obey an energy minimization principle. We show that this principle determines the parallel velocity in a unique way as the parallel gradient along the potential level sets of a velocity potential (not to be confused with the confinement potential). This velocity potential is found by inverting a Laplace-Beltrami operator on each of the level sets. 

We will show that in general, it is not possible to neglect the tangential component of the velocity. This means that the velocity at the boundary of the medium is not normal to the boundary. By contrast, the Hele-Shaw limit of the tumor models of \cite{perthame2014hele} leads to a velocity at the boundary which is normal to that boundary. Our model provides a different conclusion and consequently, brings new elements in the debate about the validity of the Darcy law, at least in its simple form when the hydraulic conductivity is a scalar. 

The medium under consideration bears analogy with a granular material. There has been considerable literature on granular media and we refer the reader to \cite{aranson2006patterns} for a review. Continuum approaches for granular media are mostly based on thermodynamical considerations (see e.g. the seminal work \cite{goodman1972continuum}). These approaches rely on the assumption that the system is at equilibrium. However, in complex media such as gels or tumors, there are momentum exchanges with the environment and energy exchanges through (bio)-chemical processes. Since these are extremely difficult to model on a first physical principle basis, we favor a heuristic approach based on the rules as described above. 

The article is structured as follows. In Section \ref{sec_framework} we summarize the main results of our work and provide a detailed discussion and directions for future work. The following sections are devoted to the proofs. The case of the mechanical equilibrium is dealt with in Section \ref{sec:minimization}. Then, the time dependent problem is investigated with first the determination of the normal velocity in Section \ref{sec:velocity} and then that of the tangential velocity in Section \ref{sec:consistency_continuity}. A short conclusion is drawn in Section \ref{sec:conclu}.

\setcounter{equation}{0}
\section{Framework,  main results and discussion}
\label{sec_framework}

\subsection{Motivation: microscopic background}
\label{subsec:micro}

In this section, we motivate our approach by proposing a model of an incompressible swelling medium at the particle level. We consider a system consisting of $N$ incompressible spherical particles of positions $x_i \in {\mathbb R}^d$, $d \geq 1$, and radii $R_i>0$, for $i=1, \ldots, N$. The radii are known but the positions are the solutions of a minimization problem. Specifically, we consider that each particle is subject to a potential energy $V(x_i, R_i)$ for a given known energy function $V(x,R)$. For simplicity, we denote by ${\mathcal X} = (x_1, \ldots, x_N)$ and ${\mathcal R} = (R_1, \ldots, R_N)$. The total energy of the system is the function 
\begin{equation}
E_{\mathcal R}({\mathcal X}) = \sum_{i=1}^N V(x_i,R_i) . 
\label{eq:energy_disc}
\end{equation}

The first problem we are interested in consists of minimizing the energy (\ref{eq:energy_disc}) over a set of admissible configurations ${\mathcal X}$ corresponding to non-overlapping spheres. Specifically, we define the admissible set by
\begin{equation}
{\mathcal A}_{\mathcal R} = \big\{ {\mathcal X} \in ({\mathbb R}^d)^N \, \, \, | \, \, \, |x_i - x_j| \geq R_i + R_j, \, \forall i,j \in \{1, \ldots, N \}, \, i \not = j \big\}. 
\label{eq:admiss_set}
\end{equation}
The minimization problem consists of finding ${\mathcal X} \in ({\mathbb R}^d)^N$ which realizes
\begin{equation}
\min_{{\mathcal X} \in {\mathcal A}_{\mathcal R}} E_{\mathcal R}({\mathcal X}) .
\label{eq:ener_min_disc}
\end{equation}
This pictures the equilibrium configuration of a granular medium made of frictionless spheres in an external potential. Introducing friction or cohesion between the grains is discarded here and will be investigated in future works. Problem (\ref{eq:ener_min_disc}) has been considered numerically in \cite{degond2017damped}. This is a non-convex problem with multiple solutions. We would like to characterize the properties of a generic solution and to this end, we will consider a continuum version of it.

The second problem we consider is the introduction of time evolution dynamics in the system. This dynamics is generated by the changes over time of the particles radii $R_i(t)$, which can increase (case of a swelling material) or decrease (case of a drying material). We also allow the potential energy $V$ to depend on time. Here we will suppose that both evolutions are given. Since, the vector of the particle radii ${\mathcal R}(t)$ changes over time, the admissible set ${\mathcal A}_{{\mathcal R}(t)}$ and the potential $V(x,t,R)$ depend on time. Consequently, solutions of (\ref{eq:ener_min_disc}) will also  depend on time. Indeed, we assume that the particles stay adiabatically at a minimum of the energy (\ref{eq:energy_disc})
and that we can extract a smooth (at least differentiable) trajectory ${\mathcal X}(t)$ among the possible solutions, at least for a small interval of time. The problem is then to find the particle velocities $v_i(t) = \frac{d x_i}{dt}$, or in other words, the vector 
\begin{equation}
{\mathcal V}(t) = (v_1(t), \ldots, v_N(t)) = \frac{d {\mathcal X}}{dt}(t). 
\label{eq:velocity_disc}
\end{equation}
Again, we discard any friction or cohesion forces between the grains which could alter the time dynamics. 

A similar problem has been investigated numerically in \cite{maury2006time}. In particular, one possible algorithm is to introduce a time discretization $t^k = k \, \Delta t$ with a time step $\Delta t >0$ and assume that ${\mathcal X}^k$ is a solution of (\ref{eq:ener_min_disc}) associated to radii ${\mathcal R}^k= {\mathcal R}(t^k)$ and potential function $V^k(x,R) = V(x,t^k,R)$. Then, time is incremented by $\Delta t$ and a new minimization problem is considered associated to radii ${\mathcal R}^{k+1}$ and potential function $V^{k+1}$. Obviously, ${\mathcal X}^k$ is not a solution of this new minimization problem. So, a new solution ${\mathcal X}^{k+1}$ is sought. To single out a unique solution among the many possible solutions of the minimization problem, we select the solution ${\mathcal X}^{k+1}$ which has the smallest distance to ${\mathcal X}^k$. In this way, a discrete configuration ${\mathcal X}^{k+1}$ is found, from which a set of discrete velocities 
\begin{equation}
{\mathcal V}^k = \frac{{\mathcal X}^{k+1} - {\mathcal X}^k}{\Delta t},
\label{eq:velocity_disc_time}
\end{equation}
is found. The selection principle above leads to the velocity ${\mathcal V}^k$ of smallest possible norm among the possible candidates.  The question is whether we can find a simple expression to determine ${\mathcal V}^k$. 

Finding a simple answer to this question seems unlikely in the discrete setting, but the problem may be easier to study at the level of a coarse-grained continuum model. So, the goal of this paper is to propose such a continuum model and to show that indeed, it is possible to determine these velocities in a unique way. We would like to stress here that it is not a goal of this paper to justify the coarse-graining procedure. Rather, we are going to postulate the problem at the continuum level as an analogue of the problem at the discrete level. The investigation of the passage from the discrete to the continuum problem will be the subject of future work (see also \cite{motsch2017short} for the coarse-graining of a related model).

\subsection{General assumptions}
\label{subsec:assumpt}

We assume a medium made of discrete entities each having finite volume and minimizing a confinement energy subject to a non-overlapping (incompressibility) constraint such as described in Section \ref{subsec:micro}. Since we are aiming at a continuum description, we do not describe each particle individually but consider their number density $n(x,t)$ and their average volume $\tau(x,t) > 0$, where $x \in {\mathbb R}^d$ is the position in a $d$-dimensional space (in practice $d=1,\, 2$ or $3$) and $t \geq 0$ is the time. The non-overlapping constraint (which, at the discrete level, was expressed by the fact that ${\mathcal X}$ must belong to the admissible set ${\mathcal A}_{\mathcal R}$) is now expressed by the fact that at any given point in space and time, the volume fraction occupied by the particles $n(x,t) \tau(x,t)$ cannot exceed $1$, i.e. 
\begin{equation}
 n(x,t) \, \tau(x,t) \leq 1. 
\label{eq:nonoverlap}
\end{equation}
Thus, $\tau^{-1}(x,t)$ is the maximal allowed (packing) density of the particles. 
We assume that $\tau(x,t)$ is a given function of space and time (exactly like in the discrete setting ${\mathcal R}$ was assumed to be a function of time) and that it is defined, positive and finite irrespective of the presence of particles at $(x,t)$. The precise value of $\tau(x,t)$ in practice depends on the modelling context and will be made precise in future work. We also impose that the particle density is nonnegative:
\begin{equation}
 n(x,t) \geq 0. 
\label{eq:nonnegative}
\end{equation}

Additionally, like in the discrete case, we assume that the total number of particles $N$ is fixed, given and is constant in time, i.e. 
\begin{equation}
\int_{{\mathbb R}^d} n(x,t) \, dx = N. 
\label{eq:totmass}
\end{equation}
Again, in future work, this assumption will be removed and replaced by a model for the growth or shrinkage of the population.

\subsection{Mechanical equilibrium}
\label{subsec:equilibrium}

We are first interested by the mechanical equilibrium. Freezing the time variable $t$ for the moment, we assume that there exists a mechanical energy 
\begin{equation}
F_t[n] = \int_{{\mathbb R}^d} V(x,t,\tau(x,t)) \, n(x,t) dx, 
\label{eq:mechener}
\end{equation}
associated with a given potential $V(x,t,\tau)$, which the particles try to minimize while satisfying the non-overlapping constraint (\ref{eq:nonoverlap}), the nonnegativity constraint (\ref{eq:nonnegative}) and the total mass constraint (\ref{eq:totmass}). In other words, our goal is to solve the following minimization problem at any given time $t$: 
\begin{eqnarray}
&& \hspace{-1cm}
\mbox{Find } n(\cdot,t): \, \, x \in {\mathbb R}^d \mapsto n(x,t) \in {\mathbb R} \, \, \mbox{ a solution of: } \nonumber \\
&& \hspace{-1cm}
\min \big\{ F_t[n(\cdot,t)] \quad |\quad n(\cdot,t)\geq 0,\ \quad  n(\cdot,t) \tau \leq 1 \, \, \mbox{ and } \, \, \int_{\RR^d} n(x,t) dx =N \big\},
\label{eq_min_problem_tdepend}
\end{eqnarray}
for $\tau:(x,t)\in \R^d \times [0,\infty) \mapsto \tau(x,t)\in \R_+$ and $N>0$ given. The potential $V(x,t,\tau)$ is the continuum analog of the discrete potential $V$ of Section \ref{subsec:micro} and Eq. (\ref{eq:mechener}) is nothing but an approximation of Eq. (\ref{eq:energy_disc}) when $N$ is large, assuming that the particle positions $x_i$ are drawn randomly, independently and identically according to the probability $N^{-1} \, n(x,t) \, dx$. Obviously, whether this independence assumption holds needs to be proved but we will leave justifications of this question to future work. 

We assume that $V \geq 0$. For the simplicity of notations, we define an ``effective potential'' $W(x,t)$ by
\begin{equation}
W(x,t) = V(x,t,\tau(x,t)). 
\label{eq:def_tildeW}
\end{equation}
We assume that, for all $t \geq 0$, we have
\begin{equation}
W(x,t) \to +\infty \quad \mbox{ as } \quad |x| \to + \infty. 
\label{eq:tildeW_to_infty}
\end{equation}

In Section \ref{sec:minimization}, we will show that, under appropriate conditions on the potential $V$ including (\ref{eq:tildeW_to_infty}), the solution $n_N(x,t)$  of the minimization problem (\ref{eq_min_problem_tdepend}) (indexed by the number $N$ of particles in the system) is given by 
\begin{equation}
n_N(x,t) = \left\{ \begin{array}{lll} \displaystyle \frac{1}{\tau(x,t)}, & \mbox{ if } & x \in \Omega_N(t), \\
0, & \mbox{ if } & x \not \in \Omega_N(t),
\end{array} \right.
\label{eq:solution}
\end{equation}
where the domain $\Omega_N(t)$ is given by
\begin{equation}
\Omega_N(t) = \{ x \in {\mathbb R}^d \, \, | \, \, 0 \leq W(x,t) \leq U_N(t) \}, 
\label{eq:Omega(t)_2}
\end{equation}
and $U_N(t)$ is the unique solution of the equation
\begin{equation}
P(U_N(t),t) = N, 
\label{eq:P(U,t)=N}
\end{equation}
with $P:$ $(u,t) \in [0,\infty)^2 \mapsto P(u,t) \in [0,\infty)$ given by
\begin{equation}
P(u,t) = \int_{\{x \in {\mathbb R}^d, \, 0 \leq W(x,t) \leq u \}} \tau^{-1}(x,t) \, dx. 
\label{eq:def_P(u,t)}
\end{equation}

\begin{figure}
\centering
\includegraphics[scale=0.6]{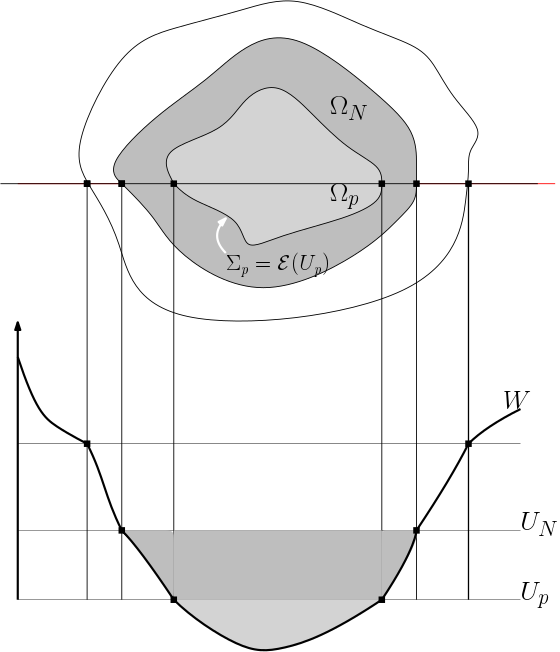}
\caption{Schematics of the filling of the potential level sets. The level set $U_N(t)$ corresponds to the filling of the potential level sets by the entire population of particles~$N$. }
\label{fig:plot1}
\end{figure}

Eq. (\ref{eq:solution}) shows, that within its support, the density saturates the congestion constraint (\ref{eq:nonoverlap}), i.e. the density is everywhere equal to the maximal allowed (packing) density $\tau^{-1}(x,t)$. Microscopically, the particles fill all the available space and it is not possible for them to increase the density any further. This is the so-called ``packing'' or ``incompressible'' state. To interpret the construction of $\Omega_N(t)$ (formulas (\ref{eq:Omega(t)_2}) to (\ref{eq:def_P(u,t)})), we introduce the level sets of the effective potential $W$. For a given $u \in [0,\infty)$ and time $t \in [0,\infty)$, the level set of $W(\cdot,t)$ corresponding to the value $u$ is defined by:
\begin{equation}
{\mathcal E}_t(u) = \{ x \in {\mathbb R}^d \, , \, \, W(x,t) = u \}. 
\label{eq:def_E_t}
\end{equation}
Eq. (\ref{eq:Omega(t)_2}) states that $\Omega_N(t)$ is bounded by the level set ${\mathcal E}_t(U_N(t))$. Formula~(\ref{eq:def_P(u,t)}) defines $P(u,t)$ as the  number of particles in the volume limited by the level $u$. Eq.~(\ref{eq:P(U,t)=N}) simply states that the level $U_N(t)$ which bounds the domain $\Omega_N(t)$ encloses the total number of particles $N$, see Fig. \ref{fig:plot1}.

Let $n_p(x,t)$ be the solution associated with a different total number of particles $p \geq 0$ with associated support $\Omega_p(t)$. 
Since $P$ is stricly increasing with respect to $u$, we have $p < N \Rightarrow U_p(t) < U_N(t)$ and so, with (\ref{eq:solution}): 
\begin{equation} 
p < N \, \Rightarrow \Omega_p(t) \subsetneq \Omega_N(t) \quad \mbox{ and } \quad n_N(\cdot,t)|_{\Omega_p(t)} = n_p(\cdot,t). 
\label{eq:induction}
\end{equation}
Additionally, 
We introduce the domain boundary $\Sigma_p(t)$ of $\Omega_p(t)$. With (\ref{eq:Omega(t)_2}) and (\ref{eq:def_E_t}), we have
\begin{equation}
\Sigma_p(t) = \partial \Omega_p(t) = \{ x \in {\mathbb R}^d \, \, | \, \, W(x,t) = U_p(t) \} = {\mathcal E}_t(U_p(t)). 
\label{eq:def_Sigma_p}
\end{equation}
This surface will play a crucial role in the definition of the dynamics below. Here, we just remark that, as a consequence of (\ref{eq:Omega(t)_2}),
\begin{equation}
\Omega_N(t) = \bigcup_{p \uparrow N} \Sigma_p(t), 
\label{eq:union_Sigma}
\end{equation}
see Fig. \ref{fig:plot1}.

\subsection{Motion under volume growth in non-swapping condition}
\label{subsec:motion} 

\begin{figure}
\centering
\includegraphics[scale=0.35]{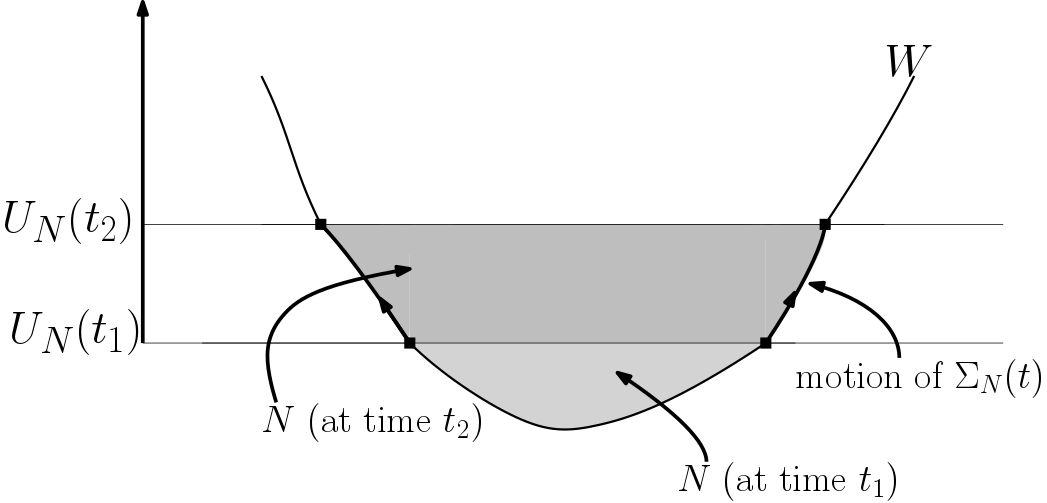}
\caption{Schematics of the motion of the medium between two instants $t_1, t_2$ where $\tau(\cdot, \tau_2) > \tau (\cdot, \tau_1)$.}
\label{fig:plot2}
\end{figure}

Now, we turn our attention towards a dynamic situation where the average volume occupied by the particles $\tau(x,t)$ at point $(x,t)$ may vary in time due to either their swelling or drying, described respectively by a time-increasing or decreasing average volume $\tau(x,t)$. We also allow for a possible time-dependence of the confinement potential function $V(x,t,\tau)$. We assume that at any given time $t$, the medium is at mechanical equilibrium as described in the previous section. So, the time variations of $\tau$ and $V$ induce an evolution of the density $n$ and of the material interface $\Omega_N(t)$ in an adiabatic way, i.e. the system follows a trajectory which is a time-continuous sequence of mechanical equilibria, see Fig. \ref{fig:plot2}. We are interested by the motion of the material-vacuum interface $\Omega_N(t)$ but also, more importantly, by the motion of the medium itself. More precisely, we would like to define a continuum velocity $v(x,t)$, $x \in \Omega_N(t)$  such that the continuity equation 
\begin{equation}
\partial_t n + \nabla \cdot (nv)=0,
\label{eq:continuity}
\end{equation}
is satisfied with the solution $n = n_N$ in the domain $\Omega_N(t)$, where $\nabla$ indicates the spatial gradient. Since within $\Omega_N(t)$, $n_N(x,t)=\tau^{-1}(x,t)$ by virtue of  (\ref{eq:solution}), Eq.  (\ref{eq:continuity}) is an equation for $v(x,t)$, namely: 
\begin{equation}
\nabla \cdot (\tau^{-1}(x,t) \, v(x,t)) = -\partial_t \tau^{-1}(x,t), \quad x \in \Omega_N(t), \quad t \geq 0. 
\label{eq:continuity_2}
\end{equation} 
However, it is a scalar equation for the vector quantity $v(x,t)$ and only fully determines~$v$ in dimension 1.
This is exactly the statement of the ``closure problem'' discussed in \cite{ambrosi2002closure}. Here our goal is to determine the velocity $v(x,t)$ fully in any dimension, by following two principles inspired by the microscopic picture, namely, (i) the non-swapping condition and (ii) the principle of smallest  displacements. Principle (i) will determine the component of $v$ normal to the family of surfaces $(\Sigma_p(t))_{p \in (0,N]}$ while Principle (ii) will determine its tangential component to these surfaces. We will investigate the consequences of Principle (i) in the present section and defer the use of Principle (ii) to the next section. 

 The non-swapping principle (Principle (i)) postulates that the level sets of the potential constrain the dynamics of the particles. More precisely, it postulates that two neighboring  particles that are on a same level set at one time will continue to be on the same level set at future times, while those on different level sets will continue to be on different level sets. This non-swapping assumption is a logical consequence of the fact that particles are at a packing state and cannot find enough free space to undertake a swapping manoeuvre in the normal direction.  In dimension $d=1$, we show that this assumption is always satisfied (given the assumptions made on the data) and consequently, the dynamics is fully determined by the continuity equation. By contrast, in dimension $d\geq2$, this assumption leads to a non-trivial condition that allows for the unique determination of the component of $v$ normal to the boundary $\Sigma_p(t)$ of $\Omega_p(t)$, for all $p \leq N$. To do so, we introduce 
\begin{equation}
\pi(x,t) = P (  \, W( x,t) \,, t ).
\label{eq:def_pi}
\end{equation}
This function gives the number of particles in the volume enclosed by the level set of the effective potential associated with its value at point $(x,t)$. By (\ref{eq:def_Sigma_p}), we have 
\begin{equation}
\Sigma_p(t) = \{x \in {\mathbb R}^d \, \, | \, \, \pi(x,t)= p \} = \pi(\cdot, t)^{-1}(\{p\}),
\label{eq:partialOmegap}
\end{equation}
so that the family $(\Sigma_p(t))_{0 \leq p \leq N}$ is nothing but the family of level sets of the function $\pi(\cdot,t)$. We assume a non-degeneracy condition: $\nabla \pi(x,t) \not = 0$, for all $(x,t) \in {\mathbb R}^d \times [0,\infty)$. 
In geometrical language, $\pi(\cdot,t)$ endows $\Omega_N(t)$ with a fiber bundle structure with base space $(0,N]$. The vector 
\begin{equation}
\nu(x,t) = \frac{\nabla \pi (x,t)}{|\nabla \pi (x,t)|}, 
\label{eq:def_nu}
\end{equation}
defines the outward unit normal to $\Sigma_p(t)$ at $x$ with $p = \pi(x,t)$. We can decompose the velocity vector $v$ as follows: 
\begin{equation}
v(x,t) = v_\perp(x,t)  + v_\parallel(x,t), \quad v_\perp (x,t) = \big( (v \cdot \nu) \, \nu \big) (x,t), \quad v_\parallel(x,t) \cdot \nu(x,t) = 0, 
\label{eq:veldecomp}
\end{equation}
for all $x \in \Omega_N(t)$, $t \in [0,\infty)$. In the sequel, $v_\perp =|v_\perp|\nu$ will be referred to as the normal velocity (with respect to the surface $\Sigma_p$ with $p = \pi(x,t)$) and $v_\parallel$ as the tangential velocity.

\begin{figure}
\centering
\includegraphics[scale=0.6]{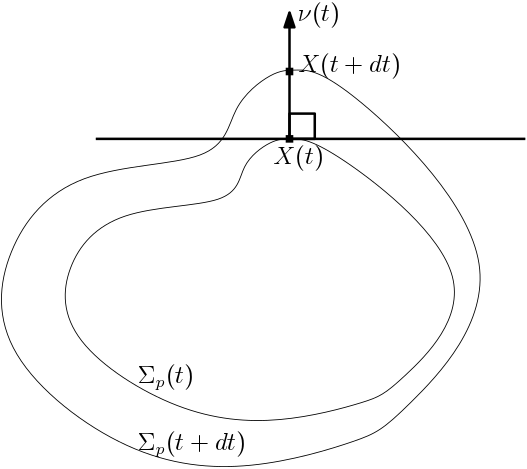}
\caption{Medium velocity in the normal direction is the velocity of $\Sigma_p$, i.e. $w_\perp~=~\frac{1}{dt}\left(X(t+dt)-X(t) \right)\cdot\nu(t)$}
\label{fig:plot3}
\end{figure}

The main consequence of the non-swapping assumption is that in dimension $d\geq 2$, it leads to the full determination of the modulus of the normal velocity $|v_\perp|=w_\perp$ as follows: 
\begin{equation}
w_\perp (x,t) = - \frac{\partial_t \pi (x,t)}{|\nabla \pi (x,t)|}, \quad x \in \Omega_N(t), \quad t \geq 0. 
\label{eq:un}
\end{equation}
This is nothing but the velocity of $\Sigma_p(t)$ in the normal direction. The interpretation is that, due to the non-swapping assumption, any particle located in the infinitesimal layer between $\Sigma_p(t)$ and $\Sigma_{p+\delta p}(t)$ with $\delta p \ll 1$ must remain in this layer and therefore, has to move with the  velocity of $\Sigma_p(t)$, see Fig. \ref{fig:plot3}.

In Section \ref{sec:velocity}, we prove that, for any velocity field satisfying (\ref{eq:un}), the left-hand side of the continuity equation (\ref{eq:continuity}) averaged on $\Sigma_p(t)$ is identically zero for any $p \leq N $ and any $t \geq 0$, namely
\begin{equation}
\big\langle \, \delta \circ \big(\pi(\cdot,t)-p\big) \, , \, \big( \partial_t n +\nabla \cdot (nv) \big)(\cdot,t) \,  \big\rangle = 0,
\label{eq:continuity_averaged}
\end{equation} 
where $\langle \cdot,\cdot \rangle$ is the duality bracket between a distribution and a smooth function. To interpret the Dirac delta in the expression above, we recall the following formula, a consequence of the so-called coarea formula:
\begin{equation} \label{eq:coarea_formula_2}
\big\langle \, \delta \circ \psi \, , \, f \,  \big\rangle = \int_{\{\psi(x) = 0 \} } f(x) \frac{dS(x)}{|\nabla \psi(x)|} , 
 \end{equation}
for any smooth functions $x \in {\mathbb R}^d \mapsto f(x), \, \psi(x) \in {\mathbb R}$, where $dS(x)$ is the euclidean surface element on the level set $\{x \in {\mathbb R}^d , \, \psi(x) = 0 \}$. The notation $(\cdot,t)$ is there to remind that the time variable $t$ is fixed when evaluating the duality bracket in (\ref{eq:continuity_averaged}). Eq. (\ref{eq:continuity_averaged}) will be an important condition for determining the tangential velocity $v_\parallel$ in the next section.

\subsection{Tangential velocity}
\label{subsec:tangential} 

To determine the tangential velocity $v_\parallel$, we apply the principle of smallest displacements (Principle (ii), see previous section). This principle suggests to determine the velocity $v_\parallel$ as the solution of a convenient energy minimization principle. It is the continuum counterpart of the principle set at the microscopic level in Section \ref{subsec:micro}, which suggested to look for the smallest velocity ${\mathcal V}^k$ among the possible ones. In the present section, we summarize the conclusions of this approach ; details and proofs can be found in Section~\ref{sec:consistency_continuity}.

First, let us make a special mention of dimension $1$, as in this case, there is no tangential velocity. So, a natural question is whether Eq. (\ref{eq:un}) is compatible with the continuity equation (\ref{eq:continuity}). In Section \ref{subsec:oneD}, we will show that this is indeed the case. This will be a consequence of (\ref{eq:continuity_averaged}).

\begin{figure}
\centering
\includegraphics[scale=0.6]{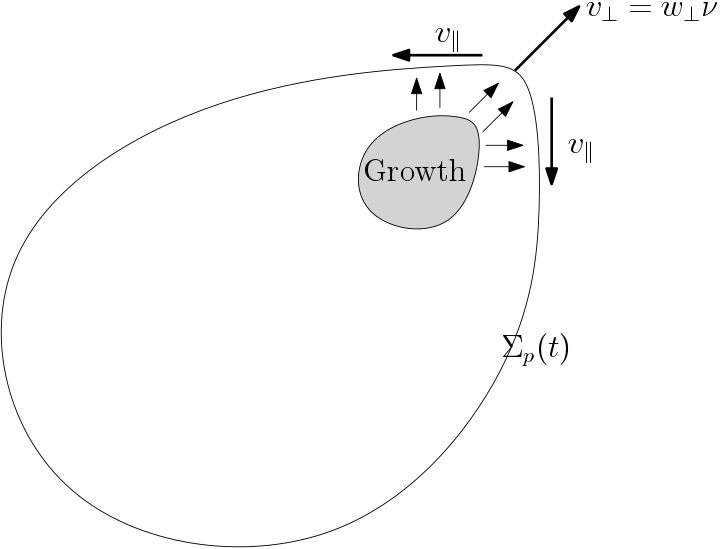}
\caption{Illustration of a need for a non-zero tangential velocity $v_\parallel$.}
\label{fig:plot4}
\end{figure}

Second, we point out that in dimension $d\geq 2$, we do need a non-zero tangential velocity $v_\parallel$ in general. Indeed, even if the choice $v = v_\perp \, \nu$ with $v_\perp$ as in (\ref{eq:un}) satisfies (\ref{eq:continuity_averaged}), it does not necessarily satisfy the continuity equation (\ref{eq:continuity}). In Section \ref{subsec:counterexample}, we will give a two-dimensional counter-example where this is indeed not true, see Fig. \ref{fig:plot4}.

So, if $d\geq 2$, in order to satisfy the continuity equation, the velocity has to incorporate a non-trivial parallel component $v_\parallel$. Using (\ref{eq:veldecomp}), the continuity equation (\ref{eq:continuity_2}) can be written
\begin{equation}
\nabla \cdot (\tau^{-1} \, v_\parallel) = f, \quad \quad f:= -\partial_t \tau^{-1} - \nabla \cdot (\tau^{-1} \, w_\perp \, \nu) , \quad \quad x \, \, \in \Omega_N(t), \, \,  t \geq 0,
\label{eq:continuity_3}
\end{equation}
and appears as a constraint on $v_\parallel$. Eq. (\ref{eq:continuity_averaged}) tells us that the average of the level sets of $f$ are all zero, namely 
\begin{equation}
\big\langle \, \delta \circ \big(\pi(\cdot,t)-p\big) \, , \, f(\cdot,t) \,  \big\rangle = 0, \quad \quad \forall (p,t) \in (0,N] \times [0,\infty). 
\label{eq:f_averaged}
\end{equation} 
In section \ref{subsec:tanvelsol}, we show that (\ref{eq:f_averaged}) is a necessary condition for the existence of a solution to (\ref{eq:continuity_3}). It is also a sufficient condition. However, in order to guarantee the uniqueness of the solution, we need to impose an additional constraint. 

Here, we add the condition that $v_\parallel$ corresponds to the minimal displacement on each of the level sets $\Sigma_p(t)$. In other words, we search for the vector fields $v_\parallel$ that minimize the parallel kinetic energy 
\begin{equation}
{\mathcal K}_{p,t}[v_\parallel] = \big\langle \, \delta \circ \big(\pi(\cdot,t)-p\big) \, ,  \, |v_\parallel (\cdot,t)|^2 \,  \big\rangle, \quad \quad \forall (p,t) \in (0,N] \times [0,\infty), 
\label{eq:kin_energy}
\end{equation} 
on all surfaces $\Sigma_p(t)$, i.e. 
\begin{equation}
v_\parallel \in \mbox{arg} \min \{ \, {\mathcal K}_{p,t}[w_\parallel] , \, \, w_\parallel \mbox{ s.t. } \nabla \cdot (\tau^{-1} \, v_\parallel) = f \, \}, \quad \quad \forall (p,t) \in (0,N] \times [0,\infty), 
\label{eq:kin_energy_min}
\end{equation} 
where we denote by arg min the set of minimizers of the expression inside the curly brackets. The expression (\ref{eq:kin_energy}) is nothing but the integral of the parallel kinetic energy density on the surface $\Sigma_p(t)$. Indeed, the parallel kinetic energy of a particle of volume $\tau$ is proportional to  $\tau |v_\parallel|^2$ but the density of such particles is proportional to $\tau^{-1}$. The contributions of the particle volume $\tau$ cancel, which leads to the expression (\ref{eq:kin_energy}). 

In section \ref{subsec:tanvelsol}, we show that such vector fields are necessarily surface gradients on the level set $\Sigma_p(t)$ of scalar functions. Specifically, we will show that (\ref{eq:kin_energy_min}) implies that there exists a scalar function $\theta(x,t)$, such that: 
\begin{equation}
v_\parallel(x,t) = - \nabla_\parallel \theta(x,t) , \quad \quad \nabla_\parallel \theta(x,t) := \nabla \theta(x,t) - \big( \nabla \theta(x,t) \cdot \nu(x,t) \big) \, \nu(x,t) , 
\label{eq:prescription_vt}
\end{equation}
where $\nabla_\parallel$ is the tangential gradient parallel to the level sets $\Sigma_p(t)$. 
With this condition, (\ref{eq:continuity_3}) becomes an elliptic equation for $\theta$ on each level set surface $\Sigma_p(t)$,  written as 
\begin{equation}
- \nabla_\parallel \cdot (\tau^{-1} \, \nabla_\parallel \theta) = f , \quad \quad  x \in \Omega_N(t), \quad t \geq 0,
\label{eq:elliptic_theta}
\end{equation}
In section \ref{subsec:tanvelsol}, this equation will be shown to have a unique solution in an appropriate function space, provided that (\ref{eq:f_averaged}) holds and that $\theta$ is sought with average zero on each level surface, namely 
\begin{equation}
\big\langle \, \delta \circ \big(\pi(\cdot,t)-p\big) \, , \, \theta(\cdot,t) \,  \big\rangle = 0, \quad \quad \forall (p,t) \in (0,N] \times [0,\infty). 
\label{eq:theta_averaged}
\end{equation} 
Indeed, (\ref{eq:elliptic_theta}) can be reformulated as the inversion of a Laplace-Beltrami operator on each of the level surfaces $\Sigma_p(t)$. Standard differential geometry (see \cite{gallot1990riemannian}, Section 4.D.2) asserts that if the solution is sought in the subspace $H^1_0(\Sigma_p(t))$ of the Sobolev space $H^1(\Sigma_p(t))$ consisting of functions satisfying the additional constraint (\ref{eq:theta_averaged}), this inversion has a unique solution. 

If the problem has spherical symmetry, i.e. if there exists ${\mathcal V}: \, (r,t,\tau) \in [0,\infty)^3 \mapsto {\mathcal V}(r,t,\tau) \in [0,\infty)$ and ${\mathcal T}: (r,t) \in [0,\infty)^2 \mapsto {\mathcal T}(r,t) \in [0,\infty)$ such that $V(x,t,\tau) = {\mathcal V} (|x|,t,\tau)$, $\tau(x,t) = {\mathcal T} (|x|,t)$, then the unique solution of (\ref{eq:elliptic_theta}), (\ref{eq:theta_averaged}) is $\theta = 0$, which shows that in this case $v_\parallel = 0$ and $v = v_\perp \nu$.

\subsection{Discussion}
\label{subsec:discussion}

First, we discuss the stationary equilibrium exposed at Section \ref{subsec:equilibrium}. The result given in (\ref{eq:solution}), (\ref{eq:Omega(t)_2}) proves that the solution of the minimization problem is unique, contrary to the discrete case exposed in Section \ref{subsec:micro}. These formulas show that the particles gradually fill the energy level sets of the effective potential $W$ by increasing values while keeping the non-overlapping condition saturated (i.e. the density being equal to the packing density). Indeed, the effective potential $W$ tends to bring all particles towards its points of global minimum. However, the non-overlapping constraint prevents the particles to pile up at these points and forces them to occupy increasingly higher potential values. They do so until the total number of particles has been exhausted. When this happens, the medium has reached its outer boundary and is therefore limited by the level set that encloses a number of particles exactly equal to the total number $N$ of available particles in the system (see Eq.~(\ref{eq:P(U,t)=N})). 

This can be compared to the process by which electrons fill energy levels in a perfect crystal at zero temperature. Electrons fill the crystal energy levels  by increasing energy due to Pauli's exclusion principle which prevents a given energy level to be occupied twice. 
The energy corresponding to the last occupied energy level is called the Fermi energy. 
The present picture is similar and $U_N(t)$ (Eq.~(\ref{eq:P(U,t)=N}))  could be viewed as the Fermi energy of our medium. The measure $\frac{dP}{du}(u,t) \, du$ (see Eq. (\ref{eq:def_P(u,t)})) which can be interpreted as the infinitesimal number of particles in a small energy interval $du$ around energy $u$ is similar to what solid-state physicists call the density-of-states, see Fig. \ref{fig:plot1}.

We now comment on the time-dependent case and the determination of the velocity in Sections \ref{subsec:motion} and \ref{subsec:tangential}. In these two sections, we provide an answer to the ``closure problem'' \cite{ambrosi2002closure}, i.e. the problem of determining the velocity field consistent with the continuity equation (\ref{eq:continuity}). This answer is different from the classical one relying on Darcy's law. Consequences of the use of Darcy's law for incompressible swelling materials can be found e.g. in \cite{perthame2014hele}. One of these is that, at the medium boundary, the continuum velocity is normal to the boundary. In the framework presented here, the velocity at the medium boundary does not have to be (and is not in general) normal to the boundary, due to the  presence of a non-trivial tangential velocity component. This discrepancy with Darcy's law may result from confinement by the external potential $V$ acting independently from the growth source modelled by $d\tau/dt$, see Fig. \ref{fig:plot4}. In \cite{perthame2014hele}, the confinement pressure is directly computed from the growth source term without any external potential $V$. Situations where confinement arises from external factors may be of importance for instance in tumor modelling when the tumor is confined by the surrounding tissue.  

\medskip
The model presented here is a building block towards a more realistic description of swelling materials such as swelling gels or tumours. This new modelling approach opens many exciting new research directions, from theory to numerics and modelling to applications. A (non-exhaustive) list of future directions which will be investigated in forthcoming works include the following: adding cell division; consider a potential $V$ that involves a contribution from particle interaction such as  attachment between nearby cells; coupling with chemical fields; introduction of boundary fuzziness; introduction of  a statistical description of particle volume sizes leading to a kinetic equation; taking into account multiple particle species; derivation from a microscopic model by coarse-graining; numerical approximation and applications to practical systems.

\smallskip
The following three sections provide the mathematical foundations of the results exposed so far. 

\setcounter{equation}{0}
\section{Equilibrium through confinement subject to volume exclusion constraint}
\label{sec:minimization}

In the present section, we provide the mathematical background to the conclusions exposed in Section \ref{subsec:equilibrium}, i.e. we determine the equilibrium configuration of the particles at a given time $t$. Throughout this section, $t$ is only a parameter, and so we will omit it in the expression of all the variables. The equilibrium configuration corresponds to minimizing the confinement energy $F[n]$ given by (\ref{eq:mechener}) subject to the volume exclusion constraint (\ref{eq:nonoverlap}), the nonnegativity constraint (\ref{eq:nonnegative}) and the total number of particles constraint (\ref{eq:totmass}). Therefore, we are led to solving the minimization problem (\ref{eq_min_problem_tdepend}) which we rewrite as follows since we omit the time-dependence:
\begin{eqnarray}
&& \hspace{-1cm}
\mbox{Find } n: \, \, x \in {\mathbb R}^d \mapsto n(x) \in {\mathbb R} \, \, \mbox{ a solution of: } \nonumber \\
&& \hspace{-1cm}
\min \big\{ F[n] \quad |\quad n\geq 0,\ \quad  n \tau \leq 1 \, \mbox{ and } \, \int_{\RR^d} n(x) dx =N \big\},
\label{eq_min_problem}
\end{eqnarray}
for $\tau:x\in \R^d \mapsto \tau(x)\in \R_+$ and $N>0$ given.
We recall the expressions (\ref{eq:def_tildeW}) of the effective potential $W$ and write $W = W(x)$ as we ignore the dependence with respect to $t$. We also recall the definition (\ref{eq:def_E_t}) of the level set of $W$ associated to the level value $u$ and we denote this level set by ${\mathcal E}(u)$, again ignoring the time-dependence. In this section we prove the following:

\begin{theorem}
\label{th:existence_minumum}
Assume the following: \\
(i) the functions $x \in {\mathbb R}^d \mapsto W(x) \in {\mathbb R}$ and $x \in {\mathbb R}^d \mapsto \tau^{-1}(x)\in {\mathbb R}$ are smooth ;  \\
(ii) $W(x)\geq 0$, $\forall x \in {\mathbb R}^d$ ; \\
(iii) $0<\tau(x)<\infty$ for all $x \in {\mathbb R}^d $ ; \\
(iv) $|\nabla W(x)|<\infty$, for all $x \in {\mathbb R}^d$ ; \\
(v) for all $u \geq 0$, the level sets ${\mathcal E}(u)$ are compact and have strictly positive $d-1$ Lebesgue surface measure; \\
(vi) $x=0$ is the only critical point of $W$ and $W(0)=0$; \\
(vii) $W(x) \to +\infty$  as $|x| \to + \infty$ ; \\
(viii) $\int_{\R^d} \tau^{-1}(t,x)\, dx > N$ for all time $t\geq 0$;\\
then, the solution of the minimization problem (\ref{eq_min_problem}) is unique and given by (\ref{eq:solution}) with the set $\Omega$ given by (\ref{eq:Omega(t)_2})-(\ref{eq:def_P(u,t)}). 
\label{prop_sol_min}
\end{theorem}

\begin{remark}
\label{rem:th_existence}
\begin{itemize}
\item[(i)] That ${\mathcal E}(u)$ is compact for all $u \geq 0$ (see Assumption (v)) follows from Assumption (vii). However, that they have strictly positive $d-1$ dimensional measure does not follow from Assumption~(vii). Conversely, Assumption (vii) does not follow from the compactness of ${\mathcal E}(u)$.

\item[(ii)] Differentiating expression (\ref{eq:def_pi}) with respect to $x$, we obtain:
\begin{equation} \label{eq:derivative_pi}
\nabla \pi(x,t)= \frac{dP}{du} (W(x,t),t) \, \nabla W(x,t).
\end{equation}
By assumption $(vi)$, $\nabla W (x,t)\neq 0$ for $x\neq 0$ and, as we will see in the proof of Th.~\ref{th:existence_minumum}, Eq. (\ref{eq:P_increasing}), it holds that 
$$\frac{dP}{du}(u)>0.$$
Therefore, from (\ref{eq:derivative_pi}) we conclude that
\begin{equation} \label{eq:derivation_nondegeneracy}
\nabla \pi (x,t) \neq 0, \quad
\mbox{for }x\neq 0,
\end{equation} 
which is a non-degeneracy condition that we will use in the sequel. Moreover, by Assumption (i), using Eq. (\ref{eq:def_pi}), we have that $\pi$ is also smooth. 

\item[(iii)] A more general form of the coarea formula (\ref{eq:coarea_formula}) (see Ref. \cite{federer2014geometric}) would allow us to extend the results with weaker assumptions than $(vi)$ or without having to assume that $\nabla \pi \neq 0$. However, to keep the presentation simple, we do not follow this path here. Indeed, with assumption $(vi)$ we ensure that $\Omega$ stays connected. If we had, say, two connected components, the global minimisation problem (\ref{eq_min_problem}) would fix the number of particles in each of the connected components, which is unrealistic, as we may expect that these two numbers could a priori be chosen independently. 
\end{itemize}
\end{remark}

Before proving Theorem \ref{prop_sol_min} we first prove the following:
\begin{lemma}
Suppose the assumptions of Theorem \ref{prop_sol_min} hold. Then, a solution $n$ of the minimization problem (\ref{eq_min_problem}) is such that, for all $x \in {\mathbb R}^d$, 
\begin{equation}
\mbox{ either } \quad n (x) \, \tau(x) = 1 \quad  \mbox{ or } \quad n(x) = 0.
\label{eq:ntau=1}
\end{equation}
\label{lem_sol_mini}
\end{lemma}

\noindent
{\bf Proof:} Let $n$ be a solution to the minimization problem. Then, there exist three Lagrange-Kuhn-Tucker multipliers (see \cite[Sec. 9.2]{clarke2013functional}) $\lambda,\mu$ and $\nu$, where $\mu \in {\mathbb R}$ and $\lambda = \lambda(x) \geq 0 $ and $\nu = \nu(x) \geq 0 $ are functions satisfying: (i) $\lambda(x) = 0$ for all $x$ such that $n(x) \, \tau(x) < 1$ ; and (ii) $\nu(x) = 0$ for all $x$ such that $n(x) > 0$ ; 
such that the Euler-Lagrange equations hold: 
$$
\int W(x) \,  \delta n (x) \, dx 
= - \int \lambda(x) \, \tau(x) \, \delta n(x) \, dx + \int \nu(x) \, \delta n(x) \, dx + \mu \int \delta n(x) \, dx  ,
$$
for all small variations $\delta n(x)$ of $n(x)$. The last term corresponds to the constraint on the total mass being equal to $N$. It follows that 
$$W(x)  =  - \lambda(x) \, \tau(x) + \nu(x)  + \mu .$$
Now, suppose that $n(x') \, \tau(x') < 1$ and $n(x') >0$ for $
x'$ in a neighbourhood ${\mathcal U}$ of a point~$x$. Then, $\lambda = 0$ and $\nu = 0$ in ${\mathcal U}$ and 
$$ W(x)    = \mu  = \mbox{Constant} , \quad \forall x \in {\mathcal U} .$$
This occurrence is ruled out by Assumption (vi) of Theorem \ref{prop_sol_min}.  
Therefore, Eq. (\ref{eq:ntau=1}) must be verified. Now, thanks to condition (viii) in Th. \ref{th:existence_minumum} this is an admissible solution, which ends the proof of the Lemma. 
\endprf

Before turning to the proof of Theorem \ref{prop_sol_min}, we recall the coarea formula in its general form (formula (\ref{eq:coarea_formula_2}) is a particular case involving the Dirac delta): 
\begin{equation} \label{eq:coarea_formula} 
\int_{{\mathbb R}^d} f(x) \,  dx  = \int_{ \psi (\R^d)} \Big( \int_{\{\psi(x) = u \} } f(x) \frac{dS_u(x)}{|\nabla \psi(x)|} \Big) \, du , 
\end{equation}
where $x \in {\mathbb R}^d \mapsto \psi(x)$, $f(x) \in {\mathbb R}$ are smooth functions and $dS_u(x)$ is the euclidean surface element on the codimension-$1$ manifold $\{\psi(x) = u \}$ and $\nabla\psi$ is nowhere zero (these assumptions can be relaxed, see \cite{federer2014geometric}).  With (\ref{eq:coarea_formula_2}), we can also write (\ref{eq:coarea_formula}) as 
\begin{equation} \label{eq:coarea_formula_3}
\int_{{\mathbb R}^d} f(x) \,  dx  = \int_{\psi({\mathbb R}^d)} \big\langle \delta \circ (\psi - u) \, , \, f \big\rangle \, du . 
\end{equation}
In particular, we have 
\begin{equation} \label{eq:coarea_formula_4}
\int_{{\mathbb R}^d} f(x) \, (g \circ \psi) (x) \, dx  = \int_{\psi({\mathbb R}^d)} \big\langle \delta \circ (\psi - u) \, , \, f \big\rangle \, g(u) \,  du , 
\end{equation}
where $g:$ $\psi(\mathbb R^d) \mapsto {\mathbb R}$ is a smooth function.

\bigskip
\noindent
{\bf Proof of Proposition \ref{prop_sol_min}.} Thanks to Lemma \ref{lem_sol_mini}, any solution of (\ref{eq_min_problem}) is of the form (\ref{eq:solution}) where the only unknown is the set $\Omega$. We denote by $\chi_\Omega$ the indicator function of the set $\Omega$ (we recall that the indicator function of a set $A$ is the function that takes the value $1$ on $A$ and the value $0$ on its complement set). 
Then, by the coarea formula (\ref{eq:coarea_formula_4}) applied with $f = \tau^{-1} \chi_\Omega$, $g(u)=u$ and $\psi = W$,  we get, since $n(x) = \tau^{-1}(x)$ on $\Omega$:
\begin{eqnarray}
F[n] &=& \int_\Omega W(x) \, \tau^{-1}(x) \, dx \nonumber \\
&=& \int_{{\mathbb R}^d} W(x) \, \tau^{-1}(x) \, \chi_\Omega (x) \, dx \nonumber 
\\
&=& \int_{0}^{+\infty} \big\langle \, \delta \circ ( W -u ) \, , \, \tau^{-1}\,  \chi_\Omega \, \big\rangle \, u \,  du . 
\label{eq:Ftn_decomposed}
\end{eqnarray}
Here the integration with respect to $u$ can be taken over $[0,\infty)$ thanks to Assumption (ii) of Theorem \ref{prop_sol_min}. We recall that, following (\ref{eq:coarea_formula_2})
$$\big\langle \, \delta \circ ( W -u ) \, , \, \tau^{-1}\, \,  \chi_\Omega \, \big\rangle = \int_{{\mathcal E}(u)} \frac{\tau^{-1}(x) \, \chi_\Omega (x)\, dS_u(x)}{|\nabla W(x)|}, 
$$ 
where $dS_u(x)$ is the euclidean surface element on ${\mathcal E}(u)$ and ${\mathcal E}(u)$ is defined at (\ref{eq:def_E_t}). Consequently, the only values of $\chi_\Omega (x)$ that enter the integral (\ref{eq:Ftn_decomposed}) for a  fixed value of $u$  are those taken on ${\mathcal E}(u)$. We claim that the minimum of $F[n]$  is reached if and only if the following is satisfied: (i) $\chi_\Omega (x)$ (which is equal to $0$ or $1$) is constant (i.e. either constantly $0$ or constantly $1$) on any level set ${\mathcal E}(u)$ for all $u \geq 0$; (ii) there exists $U>0$ such that $\chi_\Omega (x) = 1$ on ${\mathcal E}(u)$ for all $u$ such that $0 \leq u \leq U$ and $\chi_\Omega (x) = 0$ for $u \geq U$. Equivalently, these two conditions  put together mean that $\chi_\Omega (x)$ can be written:
\begin{equation}
\chi_\Omega (x) = \chi_{[0,U]} (W(x)), \quad \mbox{i.e.} \quad \chi_\Omega = \chi_{[0,U]} \circ W.
\label{eq:chi_Omega(t)}
\end{equation}
It follows that (thanks to (\ref{eq:coarea_formula_4})
\begin{eqnarray}
&&\hspace{-1cm}
\big\langle \, \delta \circ ( W -u ) \, , \,  \tau^{-1} \, \chi_\Omega \, \big\rangle 
= \chi_{[0,U]} (u) \, \, \big\langle \, \delta \circ ( W -u ) \, , \,  \tau^{-1} \, \big\rangle, \label{eq:integral_xi_Omega}
\end{eqnarray}
and 
\begin{equation} 
F[n] = \int_0^{U} \big\langle \, \delta \circ ( W-u ) \, , \,  \tau^{-1} \, \big\rangle \, \,  u \, du . 
\label{eq:express_F_t}
\end{equation}

Assuming this result for a while, i.e., that $U$ satisfying (\ref{eq:chi_Omega(t)}) exists, we show that $U$ is uniquely determined by the total number of particles constraint (\ref{eq:totmass}). Using (\ref{eq:chi_Omega(t)}) and the fact that on $\Omega$, $n(x) = \tau^{-1}(x)$, we can compute the total mass as follows:
\begin{eqnarray*} 
N &=& \int_\Omega \tau^{-1}(x) \, dx  \\
&=& \int_{{\mathbb R}^d} \tau^{-1}(x) \, \chi_\Omega (x) \, dx  \\
&=& \int_{{\mathbb R}^d} \tau^{-1}(x) \, \, \chi_{[0,U]} ( W(x)) \, dx \\
&=& \int_{\{ x \in {\mathbb R}^d \, , \, \, 0 \leq W(x) \leq U \}} \tau^{-1}(x) \, dx \\
&=& P(U),
\end{eqnarray*}
where the function $P$ (for fixed time $t$) is defined by (\ref{eq:def_P(u,t)}). This leads to Eq. (\ref{eq:P(U,t)=N}) for the determination of $U$. Note that $P(U) < \infty$ for any $U \geq 0$ by Assumption (vii). 

Now, Eq. (\ref{eq:P(U,t)=N}) has a unique solution. Indeed, using the coarea formula again, we have 
\begin{eqnarray*} 
P(u) &=& \int_0^u \big\langle \, \delta \circ ( W - u' ) \, , \,  \tau^{-1} \, \big\rangle \, du' . 
\end{eqnarray*}
Therefore, using (\ref{eq:coarea_formula_2}) and recalling the definition (\ref{eq:def_E_t}) of ${\mathcal E}(u)$, we have
\begin{eqnarray*} 
\frac{d P}{d u}(u) &=& \big\langle \, \delta \circ ( W - u ) \, , \,  \tau^{-1} \, \big\rangle \\ 
&=& \int_{{\mathcal E}(u)} \tau^{-1}(x) \, \frac{dS_u(x)}{|\nabla W(x)|}.
\end{eqnarray*}
From Assumptions (i) and (iii) to  (iv) and (vi) of Theorem (\ref{prop_sol_min}), there exists $C_u >0$ such that $\tau^{-1}(x) \, |\nabla W(x)|^{-1} \geq C_u >0$ on ${\mathcal E}(u)$. Thus, by Assumption (v) of Theorem (\ref{prop_sol_min}),
\begin{eqnarray} \label{eq:P_increasing} 
\frac{d P}{d u}(u) &\geq& C_u \int_{{\mathcal E}(u)} dS_u(x) >0 .
\end{eqnarray}
Consequently, $P$ is a strictly increasing function and there exists a unique $u=U$ such that (\ref{eq:P(U,t)=N}) holds.

We now show (\ref{eq:Omega(t)_2}). Denote by $\Omega_0$ the set defined by (\ref{eq:Omega(t)_2}) and by $n_0$ the corresponding density given by (\ref{eq:solution}). Taking $\chi_\Omega$ not of the form (\ref{eq:Omega(t)_2}), we show that the corresponding density $n$ has energy strictly larger than that of $n_0$, i.e. $F[n] > F[n_0]$. This incidentally shows the uniqueness of the solution of the minimization problem as from Lemma \ref{lem_sol_mini}, it must be of the form (\ref{eq:solution}) for some set $\Omega$ and if $\Omega \not = \Omega_0$, then, its energy is strictly larger than that obtained with $\Omega_0$.

Taking $\Omega \not = \Omega_0$ means that at least one of the subsets 
$$\omega_1 = \{x \in {\mathbb R}^d, \mbox{ such that } W(x) \leq U \mbox{ and } \chi_\Omega =0 \}, $$ 
or 
$$\omega_2 = \{x \in {\mathbb R}^d, \mbox{ such that } W(x) > U \mbox{ and } \chi_\Omega =1 \},  $$
contains a non-zero number of particles (i.e. has non-zero measure for the measure $\tau^{-1}(x) \, dx$). We now show that they both contain a non-zero number of particles and that these numbers are the same by the total number of particles constraint (\ref{eq:totmass}). 
Indeed, we note that  
\begin{eqnarray}
&&\hspace{-1cm}
\Omega_0 \setminus \omega_1 = \Omega \setminus \omega_2 = \{x \in {\mathbb R}^d, \mbox{ such that } W(x) \leq U \mbox{ and } \chi_\Omega=1 \}. 
\label{eq:tildeom}
\end{eqnarray}
Denote this set by $\tilde \omega$. Then, by the constraint (\ref{eq:totmass}), we can write:
$$ N = \int_{\Omega_0} \tau^{-1}(x) \, dx = \int_{\Omega} \tau^{-1}(x) \, dx.$$
Decomposing the first integral on $\omega_1$ and $\tilde \omega$ (which form a partition of $\Omega_0$ by (\ref{eq:tildeom})) and the second one on $\omega_2$ and $\tilde \omega$ (which similarly form a partition of $\Omega$), we get: 
$$ \int_{\omega_1} \tau^{-1}(x) \, dx + \int_{\tilde \omega} \tau^{-1}(x) \, dx = \int_{\omega_2} \tau^{-1}(x) \, dx + \int_{\tilde \omega} \tau^{-1}(x) \, dx, $$
and consequently 
\begin{equation} 
\int_{\omega_1} \tau^{-1}(x) \, dx  = \int_{\omega_2} \tau^{-1}(x) \, dx,  
\label{eq:partnumbsame_1}
\end{equation}
showing that the number of particles contained in $\omega_1$ and $\omega_2$ are the same. Note that, by the coarea formula (\ref{eq:coarea_formula_3}), we can re-write (\ref{eq:partnumbsame_1}) according to:
\begin{eqnarray}
&& \hspace{-1.5cm}
\int_{0}^{+\infty} \big\langle \, \delta \circ ( W -u ) \, , \,  ( \chi_{\omega_2} - \chi_{\omega_1} ) \, \tau^{-1} \, \big\rangle \, du = 0. 
\label{eq:partnumbsame_2}
\end{eqnarray}

Now, we have, thanks to (\ref{eq:Ftn_decomposed})
\begin{eqnarray}
&& \hspace{-1.5cm}
F[n] - F[n_0] =  \int_{0}^{+\infty} \big\langle \, \delta \circ ( W -u ) \, , \,  ( \chi_\Omega - \chi_{\Omega_0} ) \, \tau^{-1} \, \big\rangle \, u \,  du. 
\label{eq:enerdiff_1}
\end{eqnarray}
We note that
$$ \chi_{\Omega_0} = \chi_{\omega_1} + \chi_{\tilde \omega}, \quad \chi_\Omega = \chi_{\omega_2} + \chi_{\tilde \omega}. $$
So, (\ref{eq:enerdiff_1}) is written
\begin{eqnarray}
&& \hspace{-1.5cm}
F[n] - F[n_0] =  \int_{0}^{+\infty} \big\langle \, \delta \circ ( W -u ) \, , \,  ( \chi_{\omega_2} - \chi_{\omega_1} ) \, \tau^{-1} \, \big\rangle \, u \,  du. 
\label{eq:enerdiff_2}
\end{eqnarray}
But we have 
$$\omega_2 \subset \{ x \in {\mathbb R}^d \, , \, \, W(x) > U \}, \quad  \omega_1 \subset \{ x \in {\mathbb R}^d \, , \, \, W(x) \leq U \} . $$
So, we can write 
\begin{eqnarray}
\hspace{-1cm}
\int_{0}^{+\infty} \big\langle \, \delta \circ ( W -u ) \, , \,  \chi_{\omega_2} \, \tau^{-1} \, \big\rangle \, u \,  du &=& 
\int_{U}^{+\infty} \big\langle \, \delta \circ ( W -u ) \, , \,  \chi_{\omega_2} \, \tau^{-1} \, \big\rangle \, u \, du \nonumber \\
&>& U \,  \int_{U}^{+\infty} \big\langle \, \delta \circ ( W -u ) \, , \,  \chi_{\omega_2} \, \tau^{-1} \, \big\rangle  \, du \nonumber \\
&=& \, U \, \int_{0}^{+\infty} \big\langle \, \delta \circ ( W -u ) \, , \,  \chi_{\omega_2} \, \tau^{-1}\, \big\rangle \,  du, 
\label{eq:enerdiff_3}
\end{eqnarray}
and similarly

\begin{eqnarray}
\hspace{-1cm}
\int_{0}^{+\infty} \big\langle \, \delta \circ ( W -u ) \, , \,  \chi_{\omega_1} \, \tau^{-1} \, \big\rangle \, u \,  du &=& 
\int_{0}^{U} \big\langle \, \delta \circ ( W -u ) \, , \,  \chi_{\omega_1} \, \tau^{-1} \, \big\rangle \, u \, du \nonumber \\
&\leq& U \,  \int_{0}^{U} \big\langle \, \delta \circ ( W -u ) \, , \,  \chi_{\omega_1} \, \tau^{-1} \, \big\rangle  \, du \nonumber \\
&=& \, U \, \int_{0}^{+\infty} \big\langle \, \delta \circ ( W -u ) \, , \,  \chi_{\omega_1} \, \tau^{-1}\, \big\rangle \,  du, 
\label{eq:enerdiff_4}
\end{eqnarray}
Therefore,
\begin{eqnarray}
&& \hspace{-1.5cm}
F[n] - F[n_0] > U \, \int_{0}^{+\infty} \big\langle \, \delta \circ ( W -u ) \, , \,  ( \chi_{\omega_2} - \chi_{\omega_1} ) \, \tau^{-1}\, \big\rangle \, du. 
\label{eq:enerdiff_5}
\end{eqnarray}
But the integral at the right-hand side of (\ref{eq:enerdiff_5}) is equal to zero by (\ref{eq:partnumbsame_2}). Consequently, we get 
$$ F[n] > F[n_0], $$
which is the result to be proved. Note that the proof relies on the fact that the inequality in (\ref{eq:enerdiff_3}) is strict. This is only true if the support of the function 
$$ u \mapsto \big\langle \, \delta \circ ( W -u ) \, , \, \chi_{\omega_2} \tau^{-1} \,  \, \big\rangle, $$
is not reduced to $\{U\}$. But if this is the case, since the involved function is smooth, this means that it is identically equal to zero. This implies that 
$$ \int_0^\infty \big\langle \, \delta \circ ( W -u ) \, , \,  \chi_{\omega_2} \, \tau^{-1}\, \big\rangle \, du = 0, $$
and this is the total number of particles in $\omega_2$. But if there are no particles contained in $\omega_2$, that means that all particles are contained in $\Omega_0$ and therefore $\Omega = \Omega_0$. So, as soon as $\Omega \not = \Omega_0$, we have a strict inequality in (\ref{eq:enerdiff_3}). This ends the proof of Prop. \ref{prop_sol_min}. \endprf

\begin{remark}
The interpretation of (\ref{eq:express_F_t}) is as follows. The measure
$$dN(u) := \frac{dP}{du} (u) \, du = \big\langle \, \delta \circ ( W -u ) \, , \,  \tau^{-1} \, \big\rangle \, du,  $$
is the number of particles comprised between the level sets ${\mathcal E}(u)$ and ${\mathcal E}(u + du)$ (similar to the density-of-states in solid-state physics, see Section \ref{sec_framework}. In this layer, the effective potential has value $u$. So, (\ref{eq:express_F_t}) expresses that we get the total energy by summing the values of the effective potential $u$  associated to the level set ${\mathcal E}(u)$ between $0$ and $U$, weighted by the number density of particles in this level set. 
\end{remark}

\setcounter{equation}{0}
\section{Continuum velocity under non-swapping constraint}
\label{sec:velocity}

In this section, we turn our attention to a time-dynamic situation, and provide the mathematical framework to the results described in Section \ref{subsec:motion}. We consider that the average volume $\tau$ of the underlying particles in the continuum medium as well as the potential function $V$ may evolve in time. However, we assume that, during this evolution, the medium stays at mechanical equilibrium under the antagonist influences of congestion and the volume exclusion constraint at any time. Due to the time-variation of $\tau$ and $V$ the particle density $n(x,t)$ will change and we are interested in finding the velocity field $v(x,t)$ of this continuum medium. Such velocity must satisfy the continuity equation (\ref{eq:continuity}). However, this equation is a scalar equation and can only determine the vector quantity $v$ in dimension one. In dimension more than $2$, we need additional physical assumptions to determine $v$. Here, we examine what additional information on $v$ we can get from assuming that the underlying particles cannot swap their positions. We refer to Section \ref{subsec:motion} for a justification of the non-swapping assumption. 

In this section, by contrast to the previous one, we restore the time-dependence of all the quantities involved in the minimization of the mechanical energy (\ref{eq:mechener}) subject to the constraints (\ref{eq:nonoverlap}), (\ref{eq:nonnegative}), (\ref{eq:totmass}). 
We recall that under the assumptions of Theorem \ref{prop_sol_min}, the particle density $n_N(x,t)$, the unique solution of this constrained minimization problem, is given by (\ref{eq:solution}), where the domain $\Omega_N(t)$ is given by (\ref{eq:Omega(t)_2})-(\ref{eq:def_P(u,t)}). 
We also recall that in dimension $d\geq 2$, the constraint that the particles cannot swap their positions implies that those contained in the layer between two neighbouring level sets $\Sigma_p(t)$ and  $\Sigma_{p+\delta p}(t)$ with $\delta p \ll 1$ at time $t$ will remain in this layer at all times. 
Such particles must move with the layer, i.e. their normal velocity to the layer must be that of the layer or, in other words, that of the boundary $\Sigma_p(t)$.

To express this velocity, we recall the expression (\ref{eq:def_pi}) of the function $\pi(x,t)$ such that $x \in \Sigma_{\pi(x,t)}(t)$. The function $\pi(x,t)$ is the number of particles in the volume enclosed by the level set of the effective potential $W$ associated with the level value $W(x,t)$. By Eq. (\ref{eq:partialOmegap}) we also have that $\Sigma_p(t)$ is the level set of the function $\pi(\cdot,t)$. We assume the non-degeneracy condition:  
\begin{equation}
\nabla \pi(x,t) \not = 0, \quad \forall (x,t) \in {\mathbb R}^d\backslash\{0\} \times [0,\infty), 
\label{eq:nondegeneracy}
\end{equation}
which is implied by the assumptions of Th. \ref{th:existence_minumum}, see Rem. \ref{rem:th_existence} point $(ii)$.
The outward unit normal to $\Omega_p(t)$ at $x$ with $p = \pi(x,t)$ is the vector $\nu(x,t)$ defined by (\ref{eq:def_nu}) and we decompose the velocity vector $v$ according to its normal and tangential components to $\Omega_p(t)$ as defined by (\ref{eq:veldecomp}). 

\smallskip
We now recall the definition of the speed of a surface (or more generally of a co-dimension $1$ manifold). 

\begin{definition}
Consider a time-dependent smooth regular domain $\Omega(t)$ and a point $x \in \partial \Omega(t)$. Then, the speed $w_\perp(x,t)$ of the surface $\partial \Omega(t)$ at $x$ is defined as follows: define $\nu(x,t)$ the outward unit normal to $\partial \Omega(t)$ at $x$. Then, for $t'$ close to $t$, the line drawn from $x$ in the direction of $\nu(x,t)$ intersects $\partial \Omega(t')$ at a  unique point $X(t')$. Then 
\begin{equation} 
 w_\perp(x,t) = \big( \frac{d}{dt'} X(t') \big)|_{t'=t} \cdot \nu(x,t). 
\label{eq:def_wn}
\end{equation}
\label{def:bndryvel}
\end{definition}

\noindent
In the case of $\Omega_p(t)$, the speed of the surface is given in the following Lemma:

\begin{lemma}
Let $\Omega(t)=\Omega_p(t)$. Then the speed of the surface $\Sigma_p(t)$  as defined in Definition \ref{def:bndryvel} is given by
\begin{equation} 
w_\perp(x,t) =  -\frac{\partial_t \pi}{|\nabla \pi|}.
\label{eq:w&nu}
\end{equation}
\label{lem:bndryvel}
\end{lemma}

\medskip
\noindent
{\bf Proof.} We can write $\pi(X(t'),t') = p$, for all $t'$ in a small neighbourhood of $t$,  with $X(t)=x$. Therefore, using (\ref{eq:def_nu}) and (\ref{eq:def_wn}):
\begin{eqnarray}
0&=& \left.\big( \frac{d}{dt'}(\pi(X(t'),t')) \big)\right|_{t'=t} \nonumber \\
&=& \partial_t \pi (x,t) + \big( \frac{d}{dt'} X(t') \big)|_{t'=t} \cdot \nabla \pi(x,t) \nonumber \\
&=& \partial_t \pi (x,t) + \Big( \big( \frac{d}{dt'} X(t') \big)|_{t'=t} \cdot \nu(x,t) \Big) \, |\nabla \pi(x,t)| \nonumber \\
&=& \partial_t \pi (x,t) + w_\perp(x,t) \, |\nabla \pi(x,t)|, 
\label{eq:vel_bndry}
\end{eqnarray}
which leads to (\ref{eq:w&nu}) and ends the proof of the Lemma. \endprf

To define the material velocity, we will need to introduce its flow:

\begin{definition}
Given a vector field $v=v(x,t)$ which we assume continuous, bounded and $\mathcal{C}^{1}$ with respect to $x$, the flow of $v$ is the unique map $\Phi_{t}^{s}:\Omega_{N}(t)\rightarrow\Omega_{N}(s)$ such that for any $x \in \Omega_{N}(t)$, the function $\eta:$ $s\mapsto \Phi_{t}^{s}(x)$ verifies 
$$\begin{cases}
\eta(t)=x,\\
\eta'(s)=v(\eta(s),s)\ \forall s\geq 0.
\end{cases}$$
\label{def:flow}
\end{definition}

We can now define the non-swapping constraint for a velocity. 
\begin{definition}
\label{def:non-swapping}
We assume that the assumptions of Theorem \ref{prop_sol_min} are satisfied. We also assume the non-degeneracy condition (\ref{eq:nondegeneracy}). The material velocity $v(x,t)$ satisfying the same assumptions as in Def. \ref{def:flow} is said to be consistent with the non-swapping constraint if and only if for all $(x,t)$ such that $x$ is a regular point of $W(t,\cdot)$, there exists a neighborhood ${\mathcal U} \times {\mathcal V} \times {\mathcal I}$ of $(x,\pi(t,x),t)$ in $ {\mathbb R}^{d} \times [0,N] \times [0,\infty)$ and a function $H_{t}:$ $(p,s) \in {\mathcal V} \times {\mathcal I} \mapsto H_t^s(p) \in {\mathbb R}$ which is continuous and $\mathcal{C}^{1}$ with respect to $s$, such that for any $s \in {\mathcal I}$ the map $p \in {\mathcal V} \mapsto H_t^s(p) \in {\mathbb R}$ is injective and such that for all $(y,s)\in \mathcal{U}\times\mathcal{I}$, 
we have 
\begin{equation}
\label{eq:LocalNonSwappingConstraint}
\pi(\Phi_{t}^{s}(y),s)=H_{t}^{s}(\pi(y,t)),
\end{equation}
\end{definition}

\begin{figure}
\centering
\includegraphics[scale=0.55]{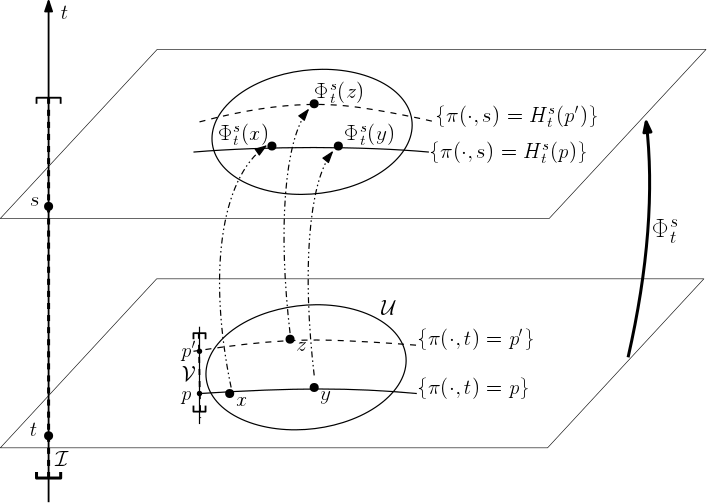}
\caption{ Schematics of the non-swapping condition in Def. \ref{def:non-swapping} }
\label{fig:plot5}
\end{figure}

\begin{remark}
Def. \ref{def:non-swapping} is illustred by Fig. \ref{fig:plot5}: Eq. (\ref{eq:LocalNonSwappingConstraint}) implies that, if at time $t$ two cells are at neighboring locations $y_1$ and $y_2$ (namely $y_1$ and $y_2$ belong to the neighborhood $\mathcal{U}$) such that they belong to the same level set, i.e. $p= \pi(y_1,t) = \pi (y_2,t)$ (respectively do not belong to the same level set i.e. $\pi(y_1,t) \not = \pi (y_2,t)$), then at time $s$ they belong to the same level set given by $\pi(\Phi_{t}^{s}(y_1),s) = \pi(\Phi_{t}^{s}(y_2),s) = H_t^s(p)$ (respectively they do not belong to the same level set i.e. $p = \pi(\Phi_{t}^{s}(y_1),s) = H_t^s(\pi(y_1,t)) \not =  p' = \pi(\Phi_{t}^{s}(y_2),s) = H_t^s(\pi(y_2,t))$ because of the injectivity of $H_t^s$). 
\end{remark}

\begin{remark}
In dimension 1 the non-swapping constraint is always satisfied and therefore carries no content. Indeed, since we suppose in Def. \ref{def:non-swapping} that $x$ is a regular point of $W(\cdot,t)$, then it is also a regular point of $\pi(\cdot,t)$, and we can locally invert $\pi(\cdot,t)_{|\mathcal{U}}:\mathcal{U}\rightarrow\mathcal{V}$. Thus we can always find a function $H_t$ satisfying Eq. (\ref{eq:LocalNonSwappingConstraint}) as
\[H_{t}^{s}(p):=\pi\left(\Phi_{t}^{s}\circ\pi_{|\mathcal{U}}(\cdot,t)^{-1}(p),s\right).\]
\label{rem:1D}
\end{remark}

Next, we give a necessary condition that the velocity $v$ has to fulfil when the evolution of $n$ is given by the continuity equation. Particularly, we show that in dimension $d\geq 2$, if a particle moves with velocity $v$ satisfying the non-swapping constraint, then the normal component of the velocity is given by the domain velocity of its level set $\Sigma_{p}(t)$ (Prop. \ref{prop:NonSwappingRelation} below), and it remains in the same level set $\Sigma_p(t)$ for all times (Prop. \ref{prop:constant_p} below). This shows that Definition \ref{def:non-swapping} ensures that a particle remains in the layer between two level sets $\Sigma_p(t)$ and $\Sigma_{p+\delta p}(t)$ at all times. More precisely:

\begin{proposition}
\label{prop:NonSwappingRelation}
Suppose that $v$ satisfies the assumptions of Def. \ref{def:flow}, verifies the non-swapping constraint as given by Def. \ref{def:non-swapping} for $d\geq 2$, and is such that
\begin{equation} \label{eq:aux_continuity}
\partial_t n + \nabla\cdot(vn)=0.
\end{equation}
Then, we have
\begin{equation}
\label{eq:NonSwappingVelocityExplicit}
v\cdot\nu=w_\perp,
\end{equation}
where $w_\perp$ is given by Eq. (\ref{eq:w&nu}) and $\nu$ by (\ref{eq:def_nu}).
\end{proposition}

\smallbreak
To prove this result, we first show the following lemma that provides a global version of the non-swapping constraint:

\begin{lemma}
\label{lemma:GlobalNonSwappingConstraint}
Let $v$ satisfy the assumptions of Def. \ref{def:flow} and verify the non-swapping constraint as expressed by Def. \ref{def:non-swapping}. We assume $d\geq 2$. Then, there exists a continuous function $h=h(p,t)$ such that for all $t,x$, 
\begin{equation}
\label{eq:GlobalNonSwappingConstraint}
(\partial_{t}+v\cdot\nabla_{x})\pi(x,t)=h(\pi(x,t),t).
\end{equation}
\label{lem:exist_h}
\end{lemma}

\medskip
\noindent
{\bf Proof.}
For all $(t,x)$ such that $x$ is not a critical point of $W(\cdot,t)$, the non-swapping constraint in Def. \ref{def:non-swapping} gives a function $H_{t}$ that verifies (\ref{eq:LocalNonSwappingConstraint}) for all $y$ in a neighbourhood $\mathcal{U}_{x}$ of $x$.
Differentiating (\ref{eq:LocalNonSwappingConstraint}) along $s$ and evaluating at $s=t$ we have:
\[(\partial_{t}+v\cdot\nabla)\pi(y,t)=\left.\partial_{s}H_{t}^{s}(\pi(y,t))\right|_{s=t},\]
for all $y\in\mathcal{U}_{x}$.
We define $h_{x}(p,t)=\left.\partial_{s}H_{t}^{s}(\pi(y,t))\right|_{s=t}$, with $\pi(y,t) = p$. We will show that this definition is independent of $x$.
Indeed, if $x,y$ are in $\Omega_{N}(t)$ (and are not critical points), and $z\in \mathcal{U}_{x}\cap\mathcal{U}_{y}$, then it must hold
\begin{equation}
\label{eq:gluingcondition}
h_{x}(\pi(z,t),t)=(\partial_{t}+v\cdot\nabla)\pi(z,t)=h_{y}(\pi(z,t),t).
\end{equation}
Now, since $d\geq 2$ and $W(\cdot,t)$ has a unique critical point (at $x=0$),  the level sets of $\pi$ are diffeomorphic to connected $(d-1)-$spheres. Using the relation (\ref{eq:gluingcondition}) and the connectivity of the level sets (since $d\geq 2$), we get that $h_{x}(p,t)=h_{y}(p,t)$ for any $x,y$ such that $h_{x}(\cdot,t)$ and $h_{y}(\cdot,t)$ are defined at $p$. Thus the functions $h_{x}$ can be glued to a single function $h=h(p,t)$ that verifies (\ref{eq:GlobalNonSwappingConstraint}). Since the functions $h_{x}$ are continuous, $h$ is continuous as well.
\endprf

\medskip
\noindent
{\bf Proof of Prop. \ref{prop:NonSwappingRelation}.} 
By lemma \ref{lemma:GlobalNonSwappingConstraint} there exists a function $h$ satisfying Eq.  (\ref{eq:GlobalNonSwappingConstraint}). This equation is equivalent to
\begin{align*}
v\cdot\nu(x,t)&=\frac{h(\pi(x,t),t)-\partial_{t}\pi(x,t)}{|\nabla\pi(x,t)|}\\
&=\frac{h(\pi(x,t),t)}{|\nabla\pi(x,t)|}+w_{\bot}(x,t),
\end{align*}
where $w_\perp$ is the normal velocity of $\Sigma_{p}(t)$ as computed in (\ref{eq:w&nu}).
Since $\int_{\Omega_{p}(t)}n\ dx=p$ by the definition of $\Omega_{p}(t)$, we deduce that
\begin{align*}
0&=\frac{d}{dt}\left(\int_{\Omega_{p}(t)}n\ dx\right)\\
&=\int_{\Omega_{p}(t)}\partial_{t}n\ dx+\int_{\Sigma_{p}(t)}n w_\perp dS(x)\\
&=\int_{\Omega_{p}(t)}-\nabla\cdot(nv)\ dx+\int_{\Sigma_{p}(t)}n w_\perp dS(x)\\
&=\int_{\Sigma_{p}(t)}n (w_\perp-v\cdot\nu) dS(x)\\
&=-h(p,t)\int_{\Sigma_{p}(t)}\frac{n}{|\nabla\pi|}dS(x).
\end{align*}
In the second line, we used the standard formula for the derivative of an integral on a time-dependent domain. The continuity equation was used in the third line, and Stokes' theorem in the fourth line. Since $n>0$ on $\Sigma_{p}(t)$, and since $\Sigma_{p}(t)$ has positive $d-1$ measure (because $p$ is not a critical value of the potential), the integral on the last line is strictly positive. We conclude that $h(p,t)=0$ for all $p>0$, and so for all $(t,x)$, we have
$$0=(\partial_t + v\cdot\nabla) \pi=\partial_t \pi + (v\cdot \nu) |\nabla \pi|,$$
which is exactly (\ref{eq:NonSwappingVelocityExplicit}) and finishes the proof. 
\endprf

As a consequence of the previous proof, we have

\begin{proposition}
Suppose $n$ satisfies the continuity equation (\ref{eq:aux_continuity}),  $v$ satisfies the non-swapping constraint as expressed in Def. \ref{def:non-swapping} and $d\geq 2$. 
Then there exists a constant $p \geq 0$, such that  $\Phi_0^t(x) \in \Sigma_p(t)$, $\forall t \geq 0$.
\label{prop:constant_p}
\end{proposition}

\medskip
\noindent
{\bf Proof.}
Let $p=\pi(x,0)$, it follows from \ref{prop:NonSwappingRelation} that
\[\frac{d}{dt}\lbrace\pi(\Phi_{0}^{t}(x),t)\rbrace=(\partial_{t}\pi+v\cdot\nabla\pi)(\Phi_{0}^{t}(x),t)=0.\]
And so $\pi(\Phi_{0}^{t}(x),t)=p$ for all $t\geq 0$, which proves the proposition. 
\endprf

\bigskip
We now show that a velocity field satisfying the non-swapping condition in dimension $d\geq 2$ (\ref{def:non-swapping}) satisfies the continuity equation averaged over all surfaces $\Sigma_p(t)$. In other words, the number of particles leaving $\Sigma_p(t)$ at a given time is exactly compensated by the number of particles arriving at $\Sigma_p(t)$. \smallbreak

\begin{theorem}
\label{prop:continuity_eq}
Under the assumptions of Theorem \ref{prop_sol_min}, let $n(x,t)$ be given by (\ref{eq:solution}). Let $v$ be a vector field such that 
\begin{equation}
v\cdot \nu = w_\perp, 
\label{eq:vnu=wn}
\end{equation} 
where $\nu$ and $w_\perp$ are given by Eqs. (\ref{eq:def_nu}) and (\ref{eq:w&nu}) respectively. Then, such vector field satisfies
\begin{equation}
\big\langle \, \delta \circ (\pi(\cdot,t) - p) \, , \,  \big( \partial_t n + \nabla \cdot (nv) \big) (\cdot,t) \, \big\rangle = 0, \quad \forall t>0, \quad \forall p \in (0, N).
\label{eq:continuity_prop}
\end{equation}
\end{theorem}

\begin{remark}
Notice that we exclude the case $p=N$ since then $n$ becomes discontinuous and the derivatives cannot be defined.
\end{remark}

\bigskip
\noindent
{\bf Proof. }  Note that we have dropped the subscript $N$ to $n_N$ for simplicity. Let $t\geq 0$, since the function \[p\mapsto \big\langle \, \delta \circ (\pi(\cdot,t) - p) \, , \,  \big( \partial_t n + \nabla \cdot (nv) \big) (\cdot,t) \, \big\rangle\] is continuous, we only need to show that for all $p\geq 0$,
\[I(p):=\int_{0}^{p}\big\langle \, \delta \circ (\pi(\cdot,t) - p') \, , \,  \big( \partial_t n + \nabla \cdot (nv) \big) (\cdot,t) \, \big\rangle dp' = 0.\]
Using Stokes' theorem, we have:
\begin{align*}
I(p) &=\int_{ \{ x \, | \, \pi(t,x)\leq p \} } \big( \partial_t n + \nabla \cdot (nv_{\bot}) \big) (x,t) \, dx\\
&=\int_{\Omega_{p}(t)} \partial_{t}n(x,t) \, dx + \int_{\Sigma_{p}(t)} nv_{\bot} \cdot \nu(x,t) \, dS(x),
\end{align*}
where $dS(x)$ is the canonical measure on $\partial\Omega(t)$. By hypothesis, $v_{\bot}\cdot\nu(x,t)$ is exactly the velocity of $\Sigma_{p}(t)$ at $(x,t)$, and so:
\begin{align*}
\int_{\Omega_{p}(t)} \partial_{t}n(x,t) \, dx + \int_{\Sigma_{p}(t)} nv_{\bot}\cdot\nu(x,t) \, dS(x)&=\frac{d}{dt}\left(\int_{\Omega_{p}(t)} n(x,t) \, dx \right)\\
&=\frac{dp}{dt}\\
&=0,
\end{align*}
where we used Eq. (\ref{eq:vel_bndry}).
So $I(p)=0$ for all $p$, which ends the proof.
\endprf

\setcounter{equation}{0}
\section{Determination of the tangential velocity}
\label{sec:consistency_continuity}

In this section, we provide the detailed mathematical discussion of the results summarized in Section \ref{subsec:tangential}.

\subsection{Dimension one}
\label{subsec:oneD}

In this section, we investigate the one-dimensional case. The non-swapping constraint is an empty constraint in this case (see Remark \ref{rem:1D}) and there is no tangential velocity. The consequence is that the dynamics of the medium is not governed by the potential (save for the determination of an integration constant), which is an important difference with the higher dimensional case. In dimension one, the continuity equation for $n$ provides a scalar differential equation for the velocity $v$, which defines it up to a constant, and this constant is determined by the boundary conditions, which indirectly involve the potential, as the following proposition shows.

\begin{proposition} 
We suppose $d=1$. Under the assumptions of Theorem \ref{prop_sol_min}, there exists a unique velocity $v$ that verifies the continuity equation (\ref{eq:aux_continuity}) and which is compatible with $n$ being a solution of the energy minimization problem, given by the conditions 
\begin{equation} \label{eq:aux_ab}
W(a(t),t)=W(b(t),t), \qquad \int_{a(t)}^{b(t)}n(x,t)dx=N,
\end{equation}
where $\Omega(t) = [a(t),b(t)]$. 
This velocity is given by 
\begin{equation}
 v(x,t)=\frac{1}{n(x,t)}\left(n(a(t),t)a'(t)-\int_{a(t)}^{x}\partial_{t}n(y,t)dy\right) ,
\label{eq:express_v_1D}
\end{equation}
where $a'(t)$ denotes the time derivative of $a(t)$ and is given by 
\begin{equation}
 a'(t) = \frac{n(b,t) \, \big( \partial_t W(b,t) - \partial_t W(a,t) \big) - \partial_x W(b,t) \int_a^b \partial_t n(x,t) \, dx}{n(b,t) \, \partial_x W(a,t) - n(a,t) \, \partial_x W(b,t)}
. 
\label{eq:express_a'}
\end{equation}
For clarity, the dependence of $a$ and $b$ on $t$ has been dropped. The expression of $b'(t)$, the time derivative of $b(t)$, is given by (\ref{eq:express_a'}) after exchanging $a$ and $b$. 
\end{proposition}
\medskip
\noindent
{\bf Proof.}
The expression of the velocity $v$  is obtained by integrating the continuity equation (\ref{eq:aux_continuity}) with respect to space on $[a(t),x]$,  noting that the velocity at $a(t)$ is precisely $a'(t)$. We just need to verify that the same property is satisfied at $b(t)$, namely that $v(b(t),t)=b'(t)$. Differentiating the second Eq. (\ref{eq:aux_ab}) with respect to $t$ gives
\begin{equation}
b'(t)n(b(t),t)-a'(t)n(a(t),t)+\int_{a(t)}^{b(t)}\partial_{t}n(y,t)dy=0.
\label{eq:deriv_part_nmbr}
\end{equation}
Using (\ref{eq:express_v_1D}), this leads to:
\[v(b(t),t)=\frac{1}{n(b(t),t)}\left(n(a(t),t)a'(t)-\int_{a(t)}^{b(t)}\partial_{t}n(y,t)dy\right)=b'(t),\]
which ends the proof. To find (\ref{eq:express_a'}) we differentiate the first Eq. (\ref{eq:aux_ab}) with respect to $t$. We find 
$$ \partial_x W(b(t),t) \, b'(t) - \partial_x W(a(t),t) \, a'(t) + \partial_t W(b(t),t) - \partial_t W(a(t),t) = 0. $$
Together with (\ref{eq:deriv_part_nmbr}), this forms a $2 \times 2$ linear system for $(a',b')$ whose solution leads to (\ref{eq:express_a'}) for $a'$ and to the corresponding expression with $a$ and $b$ exchanged for $b'$. Note that the denominator cannot be $0$ as $\partial_x W(a,t)$ and $\partial_x W(b,t)$ have opposite signs and cannot be zero as $W$ has a unique critical point which belongs to the open interval $(a(t),b(t))$. 
\endprf

\subsection{Dimension $d \geq 2$: tangential velocity is not zero in general}
\label{subsec:counterexample}

In this section, we show that in dimension $d \geq 2$ in general the velocity field must have a non-zero tangential component $v_\parallel$ to be consistent with the continuity equation. For this purpose, we provide a counter-example in dimension $d=2$ where the velocity field is defined by $v = w_\perp \nu$ with $w_\perp$ given by (\ref{eq:NonSwappingVelocityExplicit}) and which does not fulfil the continuity equation (\ref{eq:continuity}).

Indeed, consider a potential $V(x,\tau)$ which does not depend on $\tau$ and is of the form
$$V(x)= W(x) = \frac{x_2^2}{2} := \tilde W(x_2), \qquad\mbox{ for } x=(x_1,x_2)\in {\mathbb T} \times {\mathbb R},$$
and an average volume
\begin{equation} \label{eq:instance tau} 
\tau(x,t) =|x|t, \quad x\in {\mathbb T} \times {\mathbb R}, \quad t \in [0,\infty).
\end{equation}

Here ${\mathbb T} = (-1,1] \approx {\mathbb R}/{2 {\mathbb Z}}$ is the torus, i.e. we assume that all functions are $2$-periodic with respect to $x_1$ and when integrals with respect to $x_1$ are involved, they are meant over the torus ${\mathbb T}$.  
Then, by Prop. \ref{prop_sol_min} it holds that
$$n(x,t)= \frac{1}{\tau(x,t)}= \frac{1}{|x|t}, \quad x\in {\mathbb T} \times {\mathbb R}, \quad t \in [0,\infty).$$
Firstly notice that 
$$\pi(x,t)=\int_{\{ \tilde W(y_2)\leq \tilde W(x_2)\}} \tau^{-1}(y,t) \, dy :=  \tilde \pi(x_2,t),$$
so it is $x_1$-independent. The choice of $x_1$ lying in the torus ${\mathbb T}$ ensures that this integral is finite. Denoting by $(e_1,e_2)$ a cartesian basis associated to the coordinate system $(x_1,x_2)$, we get that $\nu(x,t)$ is parallel to $e_2$, i.e. 
$$\nu(x,t) = e_2 \, \mbox{ for } \, x_2 > 0, \quad \nu(x,t) = - e_2 \, \mbox{ for } \,  x_2 < 0. $$
We also have 
$$ v_\perp (x,t) = - (\partial_t \pi / |\nabla \pi| ) (x,t)  = - (\partial_t \tilde \pi / |\partial_{x_2} \tilde \pi| ) (x_2,t) : = \tilde v_\perp(x_2,t), $$ 
also only depends on $x_2$. 

This implies
$$0=\partial_t n + \nabla \cdot (n v)= \partial_t n + \partial_{x_2}(n \tilde v_\perp).$$
For the considered value of $\tau$ in (\ref{eq:instance tau}) and $x_2 > 0$, we have
$$\partial_t n + \partial_{x_2}(n \tilde v_\perp)= \frac{|x|^2 (-1+ t \, \partial_{x_2} \tilde v_\perp(x_2,t))- \tilde v_\perp(x_2) x_2 t}{|x|^3t^2}.$$
If this last expression was zero, it would imply that
$$\frac{t}{|x|^2}= \frac{-1+t \, \partial_{x_2} \tilde v_\perp(x_2,t)}{\tilde v_\perp(x_2,t) \, x_2},$$
but this cannot hold since the left-hand side depends on $x_1$ but the right-hand side  does not. Hence, we must conclude that the continuity equation is not satisfied.

\begin{remark}
The example proposed here does not satisfy the assumptions of Th. \ref{th:existence_minumum}, however it can be seen as a limiting case of $\tau^\varepsilon(x,t)=(|x|^2+\varepsilon)^{1/2} \, t$ and $V(x) = ((\varepsilon x_1^2)+x_2^2)/2$ as $\varepsilon \to 0$; and where
we have replaced assumption (vii) by periodicity conditions in the first component $x_1$.
\end{remark}

\subsection{Dimension $d \geq 2$: determination of $v_\parallel$ under principle of minimal displacement}
\label{subsec:tanvelsol}

We first show that (\ref{eq:f_averaged}) is a necessary solvability condition for (\ref{eq:continuity_3}). This is a consequence of the following lemma, in which we forget the time variable $t$: 

\begin{lemma}
Let $f$: ${\mathbb R}^d \mapsto {\mathbb R}$ be a smooth function, with $d$ a positive integer. If there exists a smooth vector field  $A$: ${\mathbb R}^d \mapsto {\mathbb R}^d$, tangent to all surfaces $\Sigma_p$, i.e.  satisfying 
\begin{equation} 
A \cdot \nabla \pi = 0, \quad \quad \mbox{in} \quad  \Omega_N,
\label{eq:Atangent}
\end{equation}
and solving the equation 
\begin{equation} 
\nabla \cdot A  = f, \quad \quad \mbox{in} \quad \Omega_N,
\label{eq:divA=f}
\end{equation}
then $f$ must be of zero-average on all level sets $\Sigma_p$, i.e. $f$ must satisfy (\ref{eq:f_averaged}). 
\label{lem:divergence}
\end{lemma}

\medskip
\noindent
{\bf Proof.} We show that if $A$: ${\mathbb R}^d \mapsto {\mathbb R}^d$ is a smooth vector field  tangent to all surfaces $\Sigma_p$, then, it satisfies 
\begin{equation}
\big\langle \, \delta \circ \big(\pi-p\big) \, , \, \nabla \cdot A  \,  \big\rangle = 0, \quad \quad \forall p \in (0,N] . 
\label{eq:averdiv}
\end{equation}
This will show the result as applying (\ref{eq:averdiv}) to (\ref{eq:divA=f}) leads to (\ref{eq:f_averaged}). To show (\ref{eq:averdiv}), we take any smooth function $g$: ${\mathbb R} \mapsto {\mathbb R}$ with compact support and compute, using (\ref{eq:coarea_formula_4}) and Green's formula: 
\begin{eqnarray*}
\int_{-\infty}^\infty g(p) \,  \big\langle \, \delta \circ \big(\pi-p\big) \, , \, \nabla \cdot A  \,  \big\rangle \, dp 
&=& \int_{{\mathbb R}^d} (g \circ \pi) (x) \, (\nabla \cdot A)(x) \, dx \\
&=& - \int_{{\mathbb R}^d} \nabla (g \circ \pi) (x) \cdot A(x) \, dx \\
&=& - \int_{{\mathbb R}^d} (g' \circ \pi) (x)  \, (\nabla \pi \cdot A)(x) \, dx \\
&=& 0,
\end{eqnarray*}
where the cancellation comes from (\ref{eq:Atangent}). This shows (\ref{eq:averdiv}) and ends the proof of the lemma. \endprf

Next, we consider the resolution of (\ref{eq:elliptic_theta}) and postpone the proof that the solution of  problem (\ref{eq:kin_energy_min}) is given by  (\ref{eq:prescription_vt}) to the end of the section. For any $(p,t) \in (0,N) \times (0,\infty)$, we note that $\Sigma_p(t) \subset \Omega_N(t)$. We denote by  ${\mathcal I}_{p,t}$: $\Sigma_p(t) \to \Omega_N(t)$ the set injection of $\Sigma_p(t)$ into $\Omega_N(t)$, i.e. for any $y \in \Sigma_p(t)$, ${\mathcal I}_{p,t} (y) = y \in \Omega_N(t)$. Now, we introduce the following change of variables. For a function $\theta$: $(x,t) \in \cup_{t \in (0,\infty)} \, \Omega_N(t) \times \{t\} \mapsto \theta(x,t) \in {\mathbb R}$, we define a function $\bar \theta$: $(p,t,y) \in \cup_{(p,t) \in (0,N) \times (0,\infty)} \,  \{ (p,t) \} \times \Sigma_p(t) \mapsto \bar \theta (p,t,y) \in (0, \infty)$ such that 
\begin{equation}
 \theta({\mathcal I}_{p,t} (y),t) = \bar \theta (p,t,y). 
\label{eq:chvar}
\end{equation}
Below, we will use that 
\begin{equation}
( \nabla_\parallel \theta) ({\mathcal I}_{p,t} (y),t) = \nabla_y \bar \theta (p,t,y), 
\label{eq:nabla_chvar}
\end{equation}
where $\nabla_y$ denote the gradient operator on the manifold $\Sigma_p(t)$. We now state the 

\begin{theorem}
Under the assumptions of Theorem \ref{prop_sol_min} and under the solvability condition (\ref{eq:f_averaged}),  Eq. (\ref{eq:elliptic_theta}) together with the zero-average constraint (\ref{eq:theta_averaged}) has a unique solution which can be written $ \theta(x,t) = \bar \theta (p,t,y)$ thanks to the change of variables (\ref{eq:chvar}), such that $\bar \theta$ belongs to the class $C^0 \big( (0,N)\times(0,\infty), H^1(\Sigma_p(t)) \big)$ where $H^1(\Sigma_p(t))$ is the Sobolev space of square integrable functions on $\Sigma_p(t)$ whose first order distributional derivatives are square integrable. 
\label{thm:existence_theta}
\end{theorem}

\medskip
\noindent
{\bf Proof.} 
Notice that $f$ (given by (\ref{eq:continuity_3})) is smooth, since $\tau^{-1}$ and $\pi$ are smooth (see Assumption (i) in Th. \ref{th:existence_minumum} and Rem. \ref{rem:th_existence} point (ii)). Taking $\psi$: $(x,t) \in {\mathbb R}^d \mapsto \psi(x,t) \in {\mathbb R}$ any smooth compactly supported function, multiplying (\ref{eq:elliptic_theta}) by $\psi$ and using Green's formula, we get:
$$\int_0^\infty \int_{{\mathbb R^d}}  \tau^{-1}(x,t) \, \nabla_\parallel \theta(x,t)  \cdot \nabla_\parallel \psi(x,t)  \, dx \, dt = \int_0^\infty \int_{{\mathbb R^d}} f(x,t) \, \psi(x,t)  \, dx \, dt , 
$$
and using  (\ref{eq:coarea_formula}), we deduce:
\begin{eqnarray}
&&\hspace{-1cm}
\int_0^\infty \int_0^N \int_{x \in \Sigma_p(t)}  \tau^{-1}(x ,t) \, \nabla_\parallel \theta(x,t)  \cdot \nabla_\parallel \psi(x,t) \, \frac{dS_{p,t}(x)}{|\nabla \pi(x,t)|} \, dp \, dt \nonumber \\
&&\hspace{3cm}
= \int_0^\infty \int_0^N \int_{x \in \Sigma_p(t)} f(x,t) \, \psi(x,t) \, \frac{dS_{p,t}(x)}{|\nabla \pi(x,t)|} \,dp \, dt , 
\label{eq:eq:theta_weak_2}
\end{eqnarray}
where $dS_{p,t}(x)$ is the euclidean surface measure on $\Sigma_p(t)$. Using the change of variable (\ref{eq:chvar}) on both $\theta$ and $\psi$, we get 
\begin{eqnarray}
&&\hspace{-1cm}
\int_0^\infty \int_0^N \int_{y \in \Sigma_p(t)}  \tau^{-1}({\mathcal I}_{p,t} (y) ,t) \, \nabla_y \bar \theta(p,t,y)  \cdot \nabla_y \bar \psi(p,t,y) \, \frac{dS_{p,t}(y)}{|\nabla \pi({\mathcal I}_{p,t} (y),t)|} \, dp \, dt \nonumber \\
&&\hspace{2cm}
= \int_0^\infty \int_0^N  \int_{y \in \Sigma_p(t)} f({\mathcal I}_{p,t} (y),t) \, \bar \psi(p,t,y) \, \frac{dS_{p,t}(y)}{|\nabla \pi({\mathcal I}_{p,t} (y),t)|} \,dp \, dt .
\label{eq:eq:theta_weak_3}
\end{eqnarray}
Since this is true for any function $\bar \psi(p,t,y)$, this implies that for any $(p,t) \in (0,N) \times (0,\infty)$, and any smooth function $\xi:$ $y \in \Sigma_p(t) \mapsto \xi(y) \in {\mathbb R}$, we have  
\begin{eqnarray}
&&\hspace{-1cm}
\int_{y \in \Sigma_p(t)}  \tau^{-1}({\mathcal I}_{p,t} (y) ,t) \, \nabla_y \bar \theta(p,t,y)  \cdot \nabla_y \xi(y) \, \frac{dS_{p,t}(y)}{|\nabla \pi({\mathcal I}_{p,t} (y),t)|} \nonumber \\
&&\hspace{3cm}
=  \int_{y \in \Sigma_p(t)} f({\mathcal I}_{p,t} (y),t) \, \xi(y) \, \frac{dS_{p,t}(y)}{|\nabla \pi({\mathcal I}_{p,t} (y),t)|} \,.
\label{eq:eq:theta_weak_6}
\end{eqnarray}
Eq. (\ref{eq:eq:theta_weak_3}) is the weak formulation of an elliptic problem posed on the closed (i.e. without boundary) smooth manifold $\Sigma_p(t)$. Reciprocally, if $y \mapsto \bar \theta(p,t,y)$ is a solution to (\ref{eq:eq:theta_weak_6}) for any $(p,t) \in (0,N) \times (0,\infty)$, then $\theta(x,t)$ constructed through (\ref{eq:chvar}) is a solution to (\ref{eq:eq:theta_weak_2}) and ultimately to (\ref{eq:elliptic_theta}). 

We now show that (\ref{eq:eq:theta_weak_6}) is equivalent to the same problem when we restrict $\xi$ to satisfy the additional constraint 
$$
\big\langle \, \delta \circ \big(\pi(\cdot,t)-p\big) \, , \, \xi \,  \big\rangle = 0, 
$$
i.e. 
\begin{equation}
\int_{y \in \Sigma_p(t)} \xi(y) \, \frac{dS_{p,t}(y)}{|\nabla \pi({\mathcal I}_{p,t} (y),t)|}  = 0 . 
\label{eq:chi_averaged}
\end{equation}
Indeed, if (\ref{eq:eq:theta_weak_6}) is satisfied for all smooth $\xi$, it is satisfied in particular for those which satisfy the additional constraint (\ref{eq:chi_averaged}). Conversely, suppose that  (\ref{eq:eq:theta_weak_6}) is satisfied for all smooth $\xi$ that satisfy (\ref{eq:chi_averaged}) and take now a smooth $\xi$ that does not satisfy (\ref{eq:chi_averaged}). We define 
$$ 
\tilde \xi (y) = \xi (y) - \frac{\int_{z \in \Sigma_p(t)} \xi(z) \, \frac{dS_{p,t}(z)}{|\nabla \pi({\mathcal I}_{p,t} (z),t)|}}{\int_{z \in \Sigma_p(t)} \, \frac{dS_{p,t}(z)}{|\nabla \pi({\mathcal I}_{p,t} (z),t)|}}.    
$$
Then, by (\ref{eq:eq:theta_weak_6}) applied with $\tilde \xi$ (which is legitimate since $\tilde \xi$ satisfies (\ref{eq:chi_averaged})), we get
\begin{eqnarray}
&&\hspace{-1cm}
\int_{y \in \Sigma_p(t)}  \tau^{-1}({\mathcal I}_{p,t} (y) ,t) \, \nabla_y \bar \theta(p,t,y)  \cdot \nabla_y \tilde \xi(y) \, \frac{dS_{p,t}(y)}{|\nabla \pi({\mathcal I}_{p,t} (y),t)|}  \nonumber \\
&&\hspace{3cm}
=  \int_{y \in \Sigma_p(t)} f({\mathcal I}_{p,t} (y),t) \, \tilde \xi(y) \, \frac{dS_{p,t}(y)}{|\nabla \pi({\mathcal I}_{p,t} (y),t)|} \,.
\label{eq:eq:theta_weak_4}
\end{eqnarray}
But since $\tilde \xi$ differs from $\xi$ by a constant on $\Sigma_p(t)$, the left-hand side of (\ref{eq:eq:theta_weak_4}) is equal to the same expression with $\xi$ instead of $\tilde \xi$. Using the assumption (\ref{eq:f_averaged}) that $f$ is of zero-average on $\Sigma_p(t)$, the right-hand side of (\ref{eq:eq:theta_weak_4}) is also equal to the same expression with $\xi$ instead of $\tilde \xi$. So, we deduce that (\ref{eq:eq:theta_weak_6}) is satisfied for all smooth $\xi$, not only those which satisfy (\ref{eq:chi_averaged}). 

So, now, we are left with solving (\ref{eq:eq:theta_weak_3}) for all smooth $\xi$ that satisfy (\ref{eq:chi_averaged}). It is time to set up functional spaces. We consider the space $L^2(\Sigma_p(t))$ of square integrable functions on $\Sigma_p(t)$ endowed with the norm 
$$ 
\|u\|_{L^2(\Sigma_p(t))} =  \Big( \int_{y \in \Sigma_p(t)} |u(y)|^2 \, \frac{dS_{p,t}(y)}{|\nabla \pi({\mathcal I}_{p,t} (y),t)|} \Big)^{1/2} \,, 
$$
and the Sobolev space $H^1(\Sigma_p(t))$ of functions $u$ of $L^2(\Sigma_p(t))$ which have first order distributional derivatives $\nabla_y u$ in $L^2(\Sigma_p(t))$, endowed with the norm 
$$ 
\|u\|_{H^1(\Sigma_p(t))} =  \Big( \|u\|^2_{L^2(\Sigma_p(t))} + \|\nabla_y u\|^2_{L^2(\Sigma_p(t))}
\Big)^{1/2} \,. 
$$
Finally, we introduce the space $H^1_0(\Sigma_p(t))$ of functions $u \in H^1(\Sigma_p(t))$ which have zero average on $\Sigma_p(t)$ i.e. such that
$$ 
\int_{y \in \Sigma_p(t)} u(y) \, \frac{dS_{p,t}(y)}{|\nabla \pi({\mathcal I}_{p,t} (y),t)|}  = 0 . 
$$
The space $H^1_0(\Sigma_p(t))$ is a closed subspace of $H^1(\Sigma_p(t))$ (because $\Sigma_p(t)$ is compact) and so, is a valid Hilbert space to apply Lax-Milgram theorem. Indeed, defining 
\begin{eqnarray*}
a(\theta, \xi )  &=& \int_{y \in \Sigma_p(t)}  \tau^{-1}({\mathcal I}_{p,t} (y) ,t) \, \nabla_y \theta(y)  \cdot \nabla_y \xi(y) \, \frac{dS_{p,t}(y)}{|\nabla \pi({\mathcal I}_{p,t} (y),t)|}, \\
\langle L, \xi \rangle &=&  \int_{y \in \Sigma_p(t)} f({\mathcal I}_{p,t} (y),t) \, \xi(y) \, \frac{dS_{p,t}(y)}{|\nabla \pi({\mathcal I}_{p,t} (y),t)|} , \\
\end{eqnarray*}
the problem of finding a solution of (\ref{eq:eq:theta_weak_3}) for all $\xi$ satisfying (\ref{eq:chi_averaged}) can be recast in the functional setting: 
\begin{eqnarray}
&&
\mbox{Find } \theta \in H^1_0(\Sigma_p(t)) \mbox{ such that } \nonumber \\
&&
\hspace{1cm} 
a(\theta, \xi ) = \langle L, \xi \rangle, \quad \forall \xi \in H^1_0(\Sigma_p(t)). 
\label{eq:LM}
\end{eqnarray}
It is clear that $a$ and $L$ are respectively a continuous bilinear form and a continuous linear form on $H^1_0(\Sigma_p(t))$. The only missing hypothesis to apply Lax-Milgram theorem is the coercivity of $a$ on $H^1_0(\Sigma_p(t))$. For this, we remark that since $\tau^{-1}$ is smooth and positive, and since $\Sigma_p(t)$ is compact, there exists $C>0$  such that $\tau^{-1}({\mathcal I}_{p,t} (y),t) \geq C >0$ for all $y \in \Sigma_p(t)$. Then, for all $\xi \in H^1_0(\Sigma_p(t))$ 
\begin{equation}
 a(\xi, \xi) \geq C \, \int_{y \in \Sigma_p(t)}  \, |\nabla_y \xi(x)|^2 \, \frac{dS_{p,t}(y)}{|\nabla \pi({\mathcal I}_{p,t} (y),t)|} := C \,  \tilde a(\xi,\xi) . 
\label{eq:achichi}
\end{equation}
The quadratic form $\tilde a(\xi,\xi)$ at the right-hand side of (\ref{eq:achichi}) is nothing but the quadratic form associated to the Laplace Beltrami operator on $\Sigma_p(t)$ endowed with the metric $g(y) = |\nabla \pi ({\mathcal I}_{p,t} (y),t)|^{-\frac{2}{d-1}} \, g_e(y)$, where $g_e(y)$ is the euclidean metric of $\Sigma_p(t)$ at point $y$. We know from the properties of the Laplace Beltrami operator on closed (i.e. without boundary) manifolds (see \cite{gallot1990riemannian}, Section 4.D.2) that its leading eigenvalue is zero, is simple and that the associated eigenfunctions are the constants. Furthermore, the eigenfunctions of this Laplace-Beltrami operator form a complete ortho-normal basis of the space $L^2(\Sigma_p(t))$. Therefore, from standard spectral theory, since $H^1_0(\Sigma_p(t))$ is the orthogonal space to the constants for the inner product of $L^2(\Sigma_p(t))$, we have
$$ 
\min_{\xi \in H^1_0(\Sigma_p(t))} \frac{\tilde a(\xi, \xi)}{ \|\xi\|_{L^2(\Sigma_p(t))}} = \lambda_1 >0,
$$ where $\lambda_1$ is the first non-zero eigenvalue of the Laplace-Beltrami operator, which is  strictly positive. Therefore, we have 
$$
a(\xi, \xi) \geq C \, \lambda_1 \|\xi\|^2_{L^2(\Sigma_p(t))}, \quad \forall \xi \in H^1_0(\Sigma_p(t)),   
$$
with $C \lambda_1 >0$, which shows the coercivity of $a$. Applying Lax-Milgram's theorem, we deduce that there exists a unique solution to (\ref{eq:LM}). Moreover, by the regularity (in $H^1$) of the solution with respect to the data, and owing to the fact that all data are smooth, we deduce that the solution $\bar \theta$ has the regularity $C^0 \big( (0,N)\times(0,\infty), H^1(\Sigma_p(t)) \big)$, which ends the proof of Theorem \ref{thm:existence_theta}. \endprf

We note that if the problem has spherical symmetry, the solution $\theta$ has also spherical symmetry, and the level sets $\Sigma_p(t)$ are spheres. Therefore, $\theta$ is constant on $\Sigma_p(t)$ but on the other hand, condition (\ref{eq:theta_averaged}) implies that its average must be zero. Therefore, the constant value of $\theta$ on $\Sigma_p(t)$ is necessarily zero. Thus, when the problem has spherical symmetry, the unique solution of (\ref{eq:elliptic_theta}), (\ref{eq:theta_averaged}) is zero, the tangential velocity $v_\parallel=0$ and the velocity $v$ is purely normal $v=w_\perp \nu$.

Now we show that the solution of minimization problem (\ref{eq:kin_energy_min}) is given by  \ref{eq:prescription_vt}). More precisely, we have the following:

\begin{proposition}
Let $v_\parallel$ be a solution of (\ref{eq:kin_energy_min}). Then, there exists a function $\theta$ such that (\ref{eq:prescription_vt}) holds. 
\label{prop:kin_energy_min}
\end{proposition}

\medskip
\noindent
{\bf Proof.} Suppose $v_\parallel= v_\parallel(x,t)$ is a solution of (\ref{eq:kin_energy_min}). Let $\delta v_\parallel = \delta v_\parallel(x,t)$ be a variation of $v_\parallel$. Then $\delta v_\parallel$ is a  tangent vector field to all level surfaces $\Sigma_p(t)$, for all $(p,t) \in (0,N) \times (0,\infty)$ and satisfies the constraint
\begin{equation}
\nabla \cdot \delta v_\parallel = 0 , \quad \forall (x,t) \in \bigcup_{t \in (0,\infty)} \, \Omega_N(t) \times \{t\}. 
\label{eq:deltavpar}
\end{equation}
Taking smooth functions $\varphi$: $(x,t) \in \cup_{t \in (0,\infty)} \, \Omega_N(t) \times \{t\} \mapsto \varphi(x,t) \in {\mathbb R}$, and $g$: $p \in (0,N) \mapsto g(p) \in {\mathbb R}$, we have, successively using Green's formula, the fact that $\delta v_\parallel$ is tangent to $\Sigma_p(t)$, and that $\nabla_\parallel (g \circ \pi) = 0$:
\begin{eqnarray*}
0 &=& \int_{\Omega_N(t)} \nabla \cdot \delta v_\parallel(x,t) \, \varphi(x,t) \, g(\pi(x,t)) \, dx  \\
&=& - \int_{\Omega_N(t)} \delta v_\parallel(x,t) \cdot \nabla (\varphi \, \, g \circ \pi)(x,t) \, dx  \\
&=& - \int_{\Omega_N(t)} \delta v_\parallel(x,t) \cdot \nabla_\parallel (\varphi \, \, g \circ \pi)(x,t) \, dx  \\
&=& - \int_{\Omega_N(t)} \delta v_\parallel(x,t) \cdot \nabla_\parallel \varphi(x,t) \, \, g (\pi(x,t)) \, dx  \\
&=& - \int_0^N \big \langle \delta \circ (\pi(\cdot ,t) - p) , \delta v_\parallel \cdot \nabla_\parallel \varphi \big \rangle \, g(p) \, dp, 
\end{eqnarray*}
where the last identity follows from (\ref{eq:coarea_formula_4}). Now, since this identity is true for all smooth functions $g(p)$, we deduce that 
\begin{eqnarray*}
0 &=&  \big \langle \delta \circ (\pi(\cdot ,t) - p) , \delta v_\parallel \cdot \nabla_\parallel \varphi \big \rangle , \quad \forall (p,t) \in (0,N) \times (0,\infty),
\end{eqnarray*}
or, using  (\ref{eq:coarea_formula_2}) and the change of variables (\ref{eq:chvar}):
\begin{eqnarray}
&&\hspace{-1.2cm} 
0=\int_{y \in \Sigma_p(t)}  \overline{\delta v}_\parallel(p,t,y)  \cdot \nabla_y \bar \varphi(p,t,y) \, \frac{dS_{p,t}( y)}{|\nabla \pi({\mathcal I}_{p,t}(y),t)|} , \, \,  \forall (p,t) \in (0,N) \times (0,\infty).
\label{eq:div_weak_1}
\end{eqnarray}

Now, the Euler-Lagrange equations of the Minimization problem (\ref{eq:kin_energy_min}) are written
\begin{eqnarray}
&& \hspace{-1cm} 
\big \langle \delta \circ (\pi(\cdot,t) - p) , v_\parallel \cdot \delta v_\parallel (\cdot, t) \big \rangle = 0, \quad \forall \, \delta v_\parallel \mbox{ tangent vector field to } \Sigma_p(t)   \nonumber \\ 
&& \hspace{2cm}
\mbox{ and satisfying (\ref{eq:deltavpar}) }, \, \, \forall (p,t) \in (0,N) \times (0,\infty),  
\label{eq:di_weak_4}
\end{eqnarray}
or, using (\ref{eq:coarea_formula_2}) and (\ref{eq:chvar}) again:
\begin{eqnarray}
0&=&\int_{y \in \Sigma_p(t)}   \bar v_\parallel (p,t,y) \cdot \overline{\delta v}_\parallel (p, t,y) \frac{dS_{p,t}(y)}{|\nabla \pi({\mathcal I}_{p,t}(y),t)|} , \quad \forall \, \delta v_\parallel \mbox{ tangent vector } \nonumber \\ 
&& \hspace{1cm}
\mbox { field to } \Sigma_p(t) \mbox{ and satisfying (\ref{eq:div_weak_1}) }, \, \, \forall (p,t) \in (0,N) \times (0,\infty),  
\label{eq:div_weak_7}
\end{eqnarray}
Eq. (\ref{eq:div_weak_7}) shows that on each surface $\Sigma_p(t)$, $\bar v_\parallel (p,t,\cdot)$ is a tangent vector field orthogonal (for the $L^2(\Sigma_p(t))$ inner product) to all tangent vector fields $\overline{\delta v}_\parallel (p, t,\cdot)$ themselves orthogonal to all gradient vector fields (by (\ref{eq:div_weak_1})). But the space of gradients of functions of $H^1(\Sigma_p(t))$ is the same as the space of gradients of functions of $H^1_0(\Sigma_p(t))$. And this latter space is closed in $L^2(\Sigma_p(t))$. This follows easily again from the coercivity of the quadratic form $\tilde a$ as proved in the proof of Theorem \ref{thm:existence_theta} (details are left to the reader). Therefore, $\bar v_\parallel (p,t,\cdot)$ being orthogonal to the orthogonal space to the gradients (and the space of gradients being closed), is itself a gradient. So, there exists a function $\bar \theta (p,t, \cdot)$ (parametrized by $(p,t) \in (0,N) \times (0,\infty)$) such that 
$$ \bar v_\parallel (p,t,y) = \nabla_y \bar \theta (p,t,y), \quad \forall y \in \Sigma_p(t), \quad \forall (p,t) \in (0,N) \times (0,\infty). $$
Defining $\theta(x,t)$ through the change of variables (\ref{eq:chvar}), we get (\ref{eq:prescription_vt}), which ends the proof. \endprf

\setcounter{equation}{0}
\section{Conclusions/perspectives}
\label{sec:conclu}

In this paper, we have proposed a new continuum model of a swelling or drying material. Two aspects have been investigated. The first one is an equilibrium problem describing particles seeking to minimize their mechanical energy subject to non-overlapping constraints. Its solution has been fully characterized. The second one is a non-equilibrium problem where we assume that the particle average volume and potential energy may vary with time and where we compute the resulting velocity applying two principles: (i)~the non swapping condition and (ii) the principle of smallest displacements. Under these two principles, the medium velocity has been fully determined. A detailed discussion has been provided and many different elaborations of the model have been proposed. In future work, we intend to progress towards the resolution of the many open problems outlined at the end of Sec. \ref{subsec:discussion}.


\bibliographystyle{abbrv}
\bibliography{biblio_tumor}

\end{document}